\documentclass[preprint, 12pt, 3p]{elsarticle_arxiv}

\usepackage{a4wide}
\usepackage{graphicx}
\usepackage{placeins}
\usepackage{subcaption}

\usepackage{amssymb}
\usepackage{amsbsy}
\usepackage{amsmath}
\usepackage{amsthm}
\usepackage{amsfonts}

\usepackage{hyperref}
\usepackage[usenames,dvipsnames]{xcolor}

\usepackage{algorithm}
\usepackage{algorithmic}

\usepackage{esint}

\journal{This manuscript version is made available under the CC-BY-NC-ND 4.0 license.}

\newcommand{\eq}[1]{(\ref{#1})}
\newcommand{\eqm}[2]{(\ref{#1})$_#2$}
\newcommand{\ol}[1]{\overline{#1}}
\newcommand{\tl}[1]{^{(#1)}}

\def\Bcal{\mathcal{B}}
\def\Ccal{\mathcal{C}}
\def\Dcal{\mathcal{D}}

\def\Gcal{\mathcal{G}}
\def\Ical{\mathcal{I}}
\def\Jcal{\mathcal{J}}
\def\Lcal{\mathcal{L}}
\def\Mcal{\mathcal{M}}
\def\Rcal{\mathcal{R}}
\def\Scal{\mathcal{S}}
\def\Qcal{\mathcal{Q}}
\def\Ucal{\mathcal{U}}
\def\Vcal{\mathcal{V}}
\def\Bcalbf{\boldsymbol{\mathcal{B}}}
\def\Dcalbf{\boldsymbol{\mathcal{D}}}
\def\Ccalbf{\boldsymbol{\mathcal{C}}}
\def\Rcalbf{\boldsymbol{\mathcal{R}}}
\def\Scalbf{\boldsymbol{\mathcal{S}}}
\def\Qcalbf{\boldsymbol{\mathcal{Q}}}

\def\gammabf{\boldsymbol{\gamma}}
\def\omegabf{\boldsymbol{\omega}}
\def\sigmabf{\boldsymbol{\sigma}}
\def\Pibf{\boldsymbol{\Pi}}

\def\ab{\boldsymbol{a}}
\def\bb{\boldsymbol{b}}
\def\eb{\boldsymbol{e}}
\def\fb{\boldsymbol{f}}
\def\nb{\boldsymbol{n}}
\def\ub{\boldsymbol{u}}
\def\vb{\boldsymbol{v}}
\def\wb{\boldsymbol{w}}

\def\yb{\boldsymbol{y}}
\def\Ab{\boldsymbol{A}}
\def\Bb{\boldsymbol{B}}
\def\Fb{\boldsymbol{F}}
\def\Hb{\boldsymbol{H}}
\def\Ib{\boldsymbol{I}}
\def\Kb{\boldsymbol{K}}
\def\Vb{\boldsymbol{V}}

\def\Aop{\pmb{\mathbb{A}}}
\def\Iop{\pmb{\mathbb{I}}}
\def\Dop{\pmb{\mathbb{D}}}

\def\HH{{\mathbb{H}}}
\def\RR{{\mathbb{R}}}
\def\SS{{\mathbb{S}}}
\def\ZZ{{\mathbb{Z}}}

\def\Lie#1{\text{\normalfont\itshape\textsterling}\kern-.1em{}_{#1}}

\def\pd{\partial}
\def\incdelta{\delta}
\def\dt{\incdelta t}
\def\veps{\varepsilon}
\def\Om{\Omega}
\def\Omt{{\Om(t)}}
\def\Omk{{\Om\tl{k}}}
\def\vphi{\varphi}

\def\Hdb{{\bf{H}}^1}
\def\HpbY{{\bf H}_\#^1(Y)}
\def\Hp{H_{\#}^1}
\def\Hpalpha{\Hp(Y_\alpha)}
\def\Hpzero{H_{\#0}^1(Y_3)}
\def\Vbzero{\Vb\kern-0.25em_0}

\def\avint{\fint}

\def\intY{\avint_Y}
\def\intYa{\avint_\alpha}
\def\intYm{\avint_{Y_3}}
\def\intYl{\avint_{Y_l}}

\def\sumalpha{\sum_{\alpha=1,2}}

\def\eff{{\rm eff}}
\def\effx{{}}
\def\new{{\rm new}}
\def\timek#1{\hat#1}

\def\trace{{\rm tr}}

\def\alin#1#2{a_{Y} \left(#1,#2\right)}
\def\blin#1#2#3{b_{Y_{#1}} \left(\bar\ub;#2,#3\right)}
\def\clin#1#2#3{c_{Y_{#1}} \left(\bar\ub;#2,#3\right)}
\def\dlin#1#2#3{d_{Y_{#1}} \left(#2,#3\right)}

\def\Tuf#1{{\mathcal{T}}_\veps{\left ({#1}\right )}} 
\def\Hpdb{{\bf{H}}_\#^1}
\def\Lb{{\bf{L}}}

\newcommand\dSy{\mathrm{\,dS}_{y}}

\newcommand\ie{{\it{i.e.~}}}
\newcommand\eg{{\it{e.g.~}}}

\newenvironment{nalign}{
    \begin{equation}
    \begin{aligned}
}{
    \end{aligned}
    \end{equation}
    \ignorespacesafterend
}

\newtheorem*{remark}{Remark}

\begin{document}

\begin{frontmatter}

\color{red}
\title{Homogenization of large deforming fluid-saturated
porous structures}

\author[1]{Vladim\'{\i}r Luke\v{s}\corref{cor1}}
\ead{vlukes@kme.zcu.cz}

\author[1]{Eduard Rohan}
\ead{rohan@kme.zcu.cz}

\address[1]{European Centre of Excellence, NTIS -- New Technologies for
Information Society, Faculty of Applied Sciences, University of West Bohemia,
Univerzitni\'{\i} 8, 30100 Plze\v{n}, Czech Republic}

\cortext[cor1]{Corresponding author}

\color{black}
  
\begin{abstract}

The two-scale computational homogenization method is proposed for
modelling of locally periodic fluid-saturated media subjected a to
large deformation induced by quasistatic loading.  The periodic
heterogeneities are relevant to the mesoscopic scale at which a double
porous medium constituted by hyperelastic skeleton and an
incompressible viscous fluid is featured by large contrasts in the
permeability. Within the Eulerian framework related to the current
deformed configuration, the two-scale homogenization approach is
applied to a linearized model discretized in time, being associated
with an incremental formulation.  For this, the equilibrium equation
and the mass conservation expressed in the spatial configuration are
differentiated using the material derivative with respect to a
convection velocity field. The homogenization procedure of the
linearized equations provides effective (homogenized) material
properties are computed to constitute the incremental macroscopic
problem. The coupled algorithm for the multiscale problem is
implemented using the finite element method. Illustrative 2D numerical
simulations of a poroelastic medium are presented including a simple
validation test.

\end{abstract}


\end{frontmatter}

\section{Introduction}

Poroelastic fluid saturated materials have been intensively investigated in
recent decades and there are various applications of classical models of porous
media in civil engineering, soil mechanics, biomechanics or tissue engineering.
The methods for modelling porous structures are usually based on
phenomenological approach suggested by Biot in \cite{Biot1941, Biot1955}, on
the mixture theory \cite{Bowen1976, Bowen1982, Bedford1983} or on volume
averaging methods \cite{Whitaker1999}. The drawback of the above methods may be
the fact, that these do not take into account the arrangement of the material
structure and usually the only geometry information are volume fractions of the
components of a heterogeneous medium.

Besides the mixture theory and other volume averaging methods, the
homogenization method based on the asymptotic analysis of the micromodel
equations with with oscillating material coefficients enable to obtain mass and
momentum balance equations relevant to the macroscopic scale as well as the
constitutive relations describing effective properties of porous media. The
two-scale homogenization method is based on separation of scales, the lower one
representing the microstructural arrangement by a periodic unit cell or by a
statistically representative volume element, that captures intrinsic material
properties of a heterogeneous structure. This method has been developed and
studied in two reference books written in the 80s by Bensoussan \textit{et al.}
\cite{Book-Bensoussan-Lions-Papanicolaou1978} emphasizing the mathematical
point of view and by Sanchez-Palencia \cite{Sanchez1980Book} emphasizing the
mechanical point of view. More recently, Auriault \textit{et al.}
\cite{Auriault2009Book} have written a reference book on the coupled phenomena
in heterogeneous media focused on the porous media. The original and
straightforward idea of Burridge and Keller \cite{Burridge-Keller-1982}
describing the fluid-structure interaction at the level of pores in the solid
phase enabled to derive the Biot model which, in the present paper, is used as
the ``micro-model'' featured by heterogeneities to be upscaled.

The computational homogenization approach, e.g. \cite{Yvonnet2019Book} is
consistent with the one based on the asymptotic analysis, but is usually
introduced on the virtual testing of a representative volume element (RVE) in
terms of strain modes, or other applied macroscopic fields. Solving so-called
local microscopic subproblems defined within RVE results in the characteristic
responses which are used to evaluate effective homogenized material parameters
involved in the macroscopic model.
The homogenization methods have widely been used for numerical simulations of
fiber \cite{Ramirez2019} or granular \cite{Miehe2010} composites, textile
structures \cite{Lomov2005}, building materials \cite{Zeman2008}, or for
modelling of mechanical behaviour of biological tissues, e.g.
\cite{Hollister1994, rohan-etal-jmps2012-bone,
rohan-cimrman-perfusionIJMCE2010, Lukes2010}. Under small strains and assuming
linear material behaviour, the total decoupling of the scales is possible and
the microstructure can be represented by one or only a few RVEs. Then the
numerical solutions are obtained by relatively low computational cost
algorithms. Considering finite strains or non-linear material dependencies, the
homogenization leads to a two-level finite element (FE$^2$) problem
\cite{FEYEL2003, OZDEMIR2008, Schroder2014} involving coupled micro-macro
analysis, where a huge number of microscopic problems have to be solved for
each iteration or time step on the macroscopic level. The numerical solution of
a coupled two-scale problem requires considerable computational power and, for
practical applications, it is necessary to use some model order reduction
techniques, as suggested e.g. in \cite{Michel20036937, Yvonnet2007,
Fritzen2013143, Sepe2013725, Eidel2019}.

In this paper, we present the two-scale homogenization approach to modelling of
the locally periodic fluid saturated porous medium, represented by Biot model
with large contrasts in the fluid permeability, subjected to finite strains.
Our interest in such structures is motivated by an effort to understand and be
able to mathematically model blood perfusion phenomena in soft biological
tissues at multiple scales. In this context, although the Biot model is usually
understood to describe behaviour at the macroscopic, or mesoscopic level, we
adhere the jargon of the two-scale homogenization employed to upscale the
microstructures. Therefore, we consider the Biot model to characterize
heterogeneous microstructures. Perfusion processes are closely related to
deformation in soft tissues, and therefore the appropriate coupled
deformation--diffusion model must be considered. Some soft tissues also undergo
large deformations and are characterized by a non-linear material behaviour, so
that the finite strain theory and a hyperelastic constitutive relation needs to
be incorporated into the proposed model. A simplified approach with the linear
kinematics and with deformation--dependent material coefficients expressed as
linear functions of the macroscopic response was reported in
\cite{Rohan-Lukes-AMC2015}. Homogenization of a hyperelastic structure with
embedded fluid-filled inclusions considering non-linear kinematics was treated
in \cite{Rohan2006-CaS} employing the updated Lagrangian formulations. An
alternative approach using the arbitrary Lagrangian Eulerian formulation was
published in \cite{brown_popov_efendiev_2013}. Here, we will follow the
linearization scheme proposed for the Biot model in our previous work
\cite{Rohan-Lukes-ADES2017}, however, we consider quasistatic loading only, so
that all inertia effects vanish. At the macro-level, the computational
algorithm is consistent with linearization of the residual formulation in the
Eulerian framework, such that the incremental scheme uses the updated
Lagrangian approach. The material derivative with respect to a convection
velocity field is used to differentiate the governing equations expressed in
the spatial configuration. The linearized system is subjected to the two-scale
homogenization and due to the proposed incremental scheme, the homogenized
material coefficients can be computed for given updated microscopic
configurations, as in the linear case,
cf.~\cite{rohan-cimrman-perfusionIJMCE2010}.

The article is organized as follows. After introducing some basic notations, in
Section~\ref{sec-incremental_def-dif}, we recall the governing equations for
the Biot model and we rewrite them in the residual weak form to which the
linearization based on the Eulerian formulation is applied. In
Section.~\ref{sec-homogenization}, we derive microscopic subproblems defined
within a RVE, we introduce the expressions for effective material coefficients
and formulate the homogenized macroscopic problem. Further, we present the time
stepping computational algorithm for the coupled micro-macro simulation. Then,
in Section~\ref{sec-num-simul}, we show a simple validation test and two
simulations demonstrating the features of the proposed multiscale model. The
summary and concluding remarks are in Section~\ref{sec-conlusion}.

\section{Basic notations}\label{sec-notations}

Throughout the paper, we shall adhere to the following notation. A point
position in a Cartesian frame is specified by $x=(x_1, x_2, x_3) \in \RR^3$,
where $\RR$ is the set of real numbers. The boldface notation for vectors $\ab
= (a_i)$ and second-order tensors $\bb = (b_{ij})$ is used. The second-order
identity tensor is denoted by $\Ib = (\delta_{ij})$. The fourth-order
elasticity tensor is denoted by $\Dop = (D_{ijkl})$. The superposed dot denotes
a derivative with respect to time. The gradient, divergence and Laplace
operators are denoted by $\nabla, \nabla \cdot$ and $\nabla^2$, respectively.
When these operators have a subscript referring to the space variable, it is
for indicating that the operator acts relatively to this space variable,
for instance $\nabla_y = (\pd_i^y) = (\pd/\pd y_i)$. The symbol dot `$\cdot$'
denotes the scalar product between two vectors and the symbol colon `$:$'
stands for scalar (inner) product of two second-order tensors, \eg $\Ab:\Bb =
A_{ij}B_{ij} = \trace[\Ab^T\Bb] = A_{ki}B_{kj}\delta_{ij}$, where
$\trace[\star]$ is the trace of a tensor and superscript $T$ in $\star^T$ is
the transposition operator. Operator $\otimes$ designates the tensor product
between two vectors, \eg $\ab\otimes\vb = (a_iv_j)$. Throughout the paper, $x$
denotes the global (``macroscopic'') coordinates, while the ``local''
coordinates $y$ describe positions within the representative unit cell
$Y\subset\RR^3$ which is introduce in the context of the locally periodic
structures. The normal vectors on a boundary of domains $\Om_\alpha$ (or
$Y_\alpha$) are denoted by $\nb^\alpha$, $\alpha = s,f$, to distinguish their
orientation outward to $\Om_\alpha$ (or $Y_\alpha$) when dealing with the
solid-fluid interfaces. By $\eb({\wb}) = 1/2(\nabla\wb + (\nabla\wb)^T)$ we
denote the infinitesimal strain tensor of a vector field $\wb$ (displacements,
or velocities). The following standard functional spaces are used: by
$L^2(\Om)$ we refer to square integrable functions defined in an open bounded
domain $\Om$; by $H^1(\Om)$ we mean the Sobolev space $W^{1,2}(\Om) \subset
L^2(\Om)$ formed by square integrable functions including their first
generalized derivatives. Bold notation is used to denote spaces of
vector-valued functions, e.g. $\Hdb(\Om)$; by subscript $_\#$ we refer to the
$Y$-periodic functions.

\section{Incremental deformation-diffusion problem}\label{sec-incremental_def-dif}

For a heterogeneous medium with a periodic microstructure, material parameters
are oscillating functions with the period proportional to the size of microscopic
heterogeneities. The period size can be express by a scale parameter $\veps$,
which will be introduced later in Section~\ref{sec-homogenization}. Due
to the oscillating material parameters, the field variables, e.g.\ displacements
and pressures, also depend on $\veps$. In this part, we do not emphasize this
scale dependence and we focus only on the introduction of the linearized
deformation-diffusion problem for a non-linear continuum. The homogenization
procedure will be applied to the linearized incremental form of the problem in
the next session.

\subsection{Governing equations for the Biot model}

The governing equations for fluid diffusion through a deforming incompressible porous
structure involve the Cauchy stress tensor $\sigmabf$ and the relative
perfusion velocity $\wb$, which are given by the following constitutive laws
\begin{nalign}\label{eq-ge01}
  \sigmabf &= -p \Ib + \sigmabf^\eff(\ub),\\
  \wb &= -\Kb \nabla p,
\end{nalign}
where $\ub$ is the displacement field, $p$ is the pore fluid pressure and $\Kb$
is the symmetric and positive definite hydraulic permeability tensor.
The Cauchy stress tensor consists of the strain dependent effective part
$\sigmabf^\eff$ and the part associated with the fluid pressure in pores. The
effective stress is related to a strain energy function which depends on the
deformation gradient $\Fb$, and thus it is a (non-linear) function of the
displacement field. Considering the neo-Hookean hyperelastic model,
$\sigmabf^\eff$ can be expressed as
\begin{nalign}\label{eq-ge01b}
  \sigmabf^\eff = \mu J^{-5/3}\, \mathrm{dev}(\bb).
\end{nalign}
Above, $J=\det(\Fb)$ is the relative volume change, $\bb = \Fb \Fb^T$ is the left
Cauchy--Green deformation tensor and $\mu$ is the shear modulus \cite{crisfield_1991}.

Deformations and fluid flow are driven by the equilibrium equation and the
volume conservation equation which read as
\begin{nalign}\label{eq-ge03}
  - \nabla\cdot \sigmabf = \fb,\\
  \nabla\cdot \dot\ub + \nabla\cdot \wb = 0,\
\end{nalign}
where $\fb$ stands for volume forces acting on the porous medium and
$\dot\ub$ is the skeleton (local) velocity, see also \cite{Rohan-Lukes-ADES2017}.

\subsection{Weak formulation}

We assume that the system \eq{eq-ge01}--\eq{eq-ge03} holds in an open bounded
domain $\Om$. By $\Om_0$ we denote the initial configuration which is
associated with material coordinates $X_i$, $i = 1,2,3$ and $\Om(t)$ is the
current configuration at time $t$ associated with the spatial coordinates $x_i$.
The domain boundary $\pd \Om$ is decomposed into disjoint parts as follows:
\begin{nalign}\label{eq-wf01}
  \pd \Om = \pd_u \Om \cup \pd_\sigma \Om,\quad \pd_u \Om \cap \pd_\sigma \Om = \emptyset.
\end{nalign}
This decompositions is reflected by the admissibility sets $\Vb$, $Q$ and corresponding
linear spaces $\Vbzero$, $Q_0$ employed in the weak formulation
\begin{nalign}\label{eq-wf02}
\Vb &= \{\vb|\; \vb = \ub^\pd \mbox{ on } \pd_u \Om\},\\
Q &= \{q|\; q = p^\pd \mbox{ on } \pd_p \Om \},\\
\Vbzero &= \{\vb|\; \vb = 0 \mbox{ on } \pd_u \Om\},\\
Q_0 &= \{q|\; q = 0 \mbox{ on } \pd_p \Om\}.
\end{nalign}
In the above definitions, the functional spaces are not specified explicitly,
but we assume sufficient regularity for all unknown and test functions. The
above spaces and sets are defined at a specific time $t$, so that we should use
the notation $\Vb(t), Q(t), \Vb_0(t), Q_0(t$ to emphasize the dependence on the
spatial configuration $\Om(t)$, whereby the boundary conditions are defined on
$\pd\Om(t)$ using time-dependent functions $\ub^\pd$ and $p^\pd$.

The state of the large deforming poroelastic medium, given by displacement and
pressure fields $\ub \in \Vb$ and $p \in Q$, is obtained by solving the non-linear
equation,
\begin{align}\label{eq-wf03}
  \Phi_t((\ub,p);(\vb,q)) = 0\quad \forall (\vb,q) \in \Vbzero(t)\times Q_0(t),
\end{align}
consisting of the equilibrium equation and the balance of the fluid content, so
that the residual function reads
\begin{nalign}\label{eq-wf04}
  \Phi_t((\ub,p);(\vb,0)) & =
  \int_{\Om(t)} \sigmabf:\nabla \vb - \int_{\Om(t)} \fb \cdot \vb\quad \forall \vb \in \Vbzero(t), \\
  \Phi_t((\ub,p);(0,q)) & =\int_{\Om(t)}\left (
  \nabla \cdot \dot \ub q + \Kb\nabla p \cdot \nabla q
  \right ) - \Jcal_t(q)\quad \forall q \in Q_0(t).
\end{nalign}
Above, functional $\Jcal_t$ includes the fluid mass sources and sinks. The
incremental form of the residual equation \eq{eq-wf03} can be derived using the
concept of the material derivative and the perturbation (velocity) field
$\Vcal$ defined in $\Omega(t)$, as discussed in \cite{Rohan-Lukes-ADES2017},
Section 3.

The residual equation \eq{eq-wf03} can be rewritten at time $t + \dt$ and approximated 
by its first order Taylor expansion evaluated at time $t$,
\begin{nalign}\label{eq-wf05}
  \Phi_{t + \dt}((\ub^\ast,p^\ast);(\vb,q))
  \approx\, & \Phi_t((\ub, p);(\vb,q)) + \incdelta\Phi_t((\ub, p);(\vb,q))
  \circ (\incdelta \ub, \incdelta p, \dt \Vcal, \dt),
\end{nalign}
where a new state $(\ub^\ast,p^\ast)$ at time $t + \dt$ is
computed for a given state $(\ub, p)$ at time $t$.
Expression
$\incdelta\Phi_t((\ub, p);(\vb,q)) \circ (\incdelta \ub, \incdelta p, \dt \Vcal)$
is the time increment due to the material derivative associated with convection field
$\Vcal$. The perturbed configuration is given by domain  $\Omega^\ast$,
\begin{nalign}\label{eq-wf06a}
  \Om(t+\dt) & \approx \Omega^\ast = \lbrace z \in \RR^3 \vert z = x + \dt \Vcal(x) \rbrace \equiv \Om(t) + \dt \{\Vcal \}_{\Om(t)},
\end{nalign}
and by the perturbed state $(\ub^\ast,p^\ast)$ expressed in terms of the state
increments $(\incdelta \ub, \incdelta p)$, such that
\begin{nalign}\label{eq-wf06}
  \ub^\ast = \ub(t + \dt) & \approx \ub(t) + \incdelta \ub, \quad \incdelta \ub = \dot\ub(t)\dt,\\
  p^\ast = p(t + \dt) & \approx p(t) + \incdelta p, \quad \incdelta p = \dot p(t)\dt.
\end{nalign}
In order to compute differential $\incdelta \Phi_t$, we employ the total time
derivative of $\Phi_t$, so that
\begin{nalign}\label{eq-wf07a}
  \frac{\incdelta \Phi_t}{\incdelta t} \approx \frac{{\rm d}}{{\rm d} t} \Phi_t \equiv \dot\Phi_t\;.
\end{nalign}

Differentiation of residual \eqm{eq-wf04}{1} yields
\begin{nalign}\label{eq-wf07}
  \int_\Omt \left(\nabla\Vcal \sigmabf + \Lie{\Vcal} \sigmabf \right):\nabla\vb - 
  \int_\Omt \left(\dot\fb \cdot \vb + \fb \cdot \vb \nabla \cdot \Vcal \right),
\end{nalign}
where $\Lie{\Vcal}$ is the Lie derivative of the Cauchy stress tensor which is
decomposed into its effective and volumetric parts. The Lie derivative of the
volumetric part $\sigmabf^p = -p\Ib$ is
\begin{nalign}\label{eq-wf08}
  \Lie{\Vcal}\sigmabf^p
    = p \left(\nabla \Vcal + (\nabla\Vcal)^T\right) - (p\nabla\cdot\Vcal + \dot{p})\Ib
\end{nalign}
and the derivative of the effective part $\sigmabf^\eff$ can be expressed in
terms of the tangential stiffness tensor $\Dop^\eff$ (Truesdell rate of the
effective stress) and linear velocity strain $\eb(\dot\ub)$,
\begin{nalign}\label{eq-wf09}
  \Lie{\Vcal}\sigmabf^\eff
    = \Dop^\eff \eb(\dot\ub), \qquad \mbox{ where } \quad \eb(\vb) = {\frac{1}{2}\left(\nabla \vb + (\nabla \vb)^T\right)}.
\end{nalign}
Upon differentiating of the volume conservation residual \eqm{eq-wf04}{2}, we get
\begin{nalign}\label{eq-wf10}
  \int_\Omt &\left(q \nabla\cdot \dot\ub + \Kb \nabla p \cdot \nabla q\right) \nabla \cdot\Vcal
  + \int_\Omt \left(q \nabla\cdot \ddot\ub + \Kb \nabla \dot p \cdot \nabla q\right)
  -\int_\Omt q \nabla \dot\ub \nabla \Vcal : \Ib\\
& - \int_\Omt \Kb \left(\nabla p \nabla \Vcal\right) \cdot \nabla q
  - \int_\Omt \Kb \nabla p \cdot \left(\nabla q \nabla \Vcal\right) 
  + \int_\Omt \dot\Kb \nabla p \cdot \nabla q - \dot\Jcal_t(q).
\end{nalign}

\subsection{Time discretization and incremental formulation}

\begin{figure}
  \centering
  \includegraphics[width=0.98\linewidth]{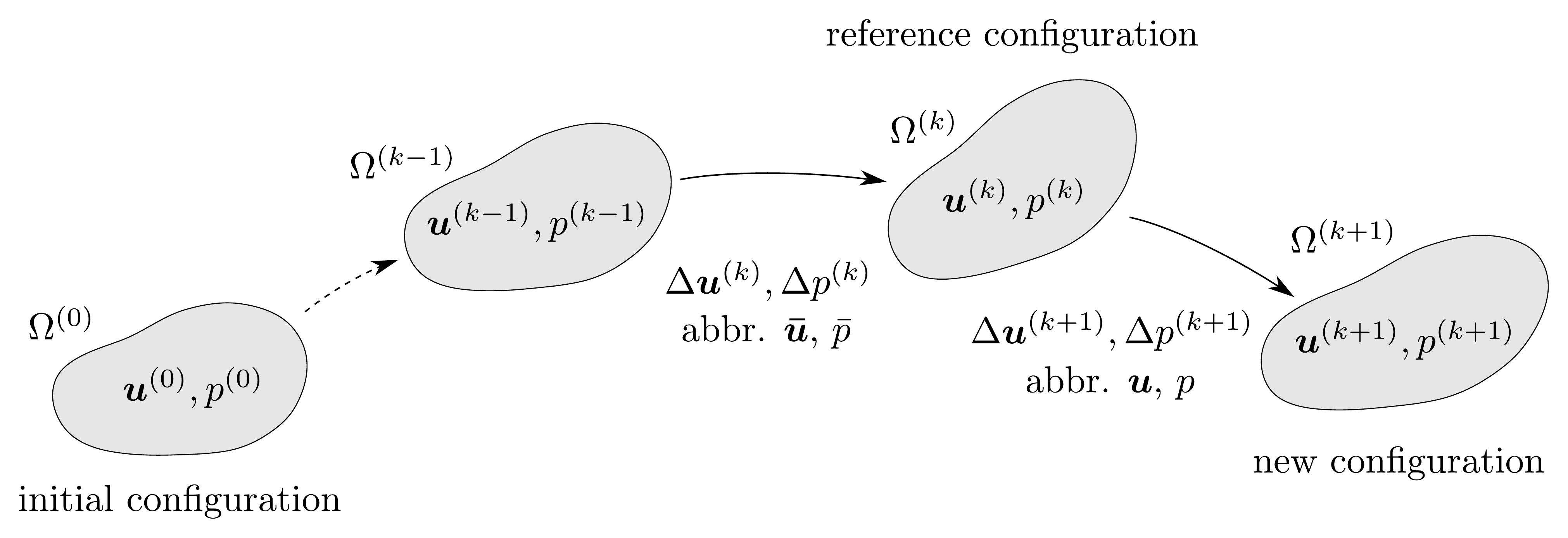}
  \caption{Updated Lagrangian Formulation: scheme of updating domain $\Omega$,
    notation of displacement and pressure fields and increments at a given time level.}\label{fig-ulf}
\end{figure}

We consider the time discretization of a time interval $]0,\bar t$ into the
time levels are $t_k = k \incdelta t$, $k = 1,2\dots$ introduced using a
constant time step. The time derivatives in \eq{eq-wf07} and \eq{eq-wf10} can
be replaced by backward differences,
\begin{nalign}\label{eq-wf11}
  \dot \ub(t_k) & = \dot \ub\tl{k} \approx (\ub\tl{k} - \ub\tl{k-1})/\incdelta t = \incdelta \ub\tl{k} / \incdelta t,\\
  \ddot \ub(t_k) & = \ddot \ub\tl{k} \approx (\incdelta \ub\tl{k} -  \incdelta \ub\tl{k-1})/(\incdelta t )^2 =
  (\ub\tl{k} - 2\ub\tl{k-1} + \ub\tl{k-2})/(\dt )^2,
\end{nalign}
where the upper index $^{(\cdot)}$ indicates the time level at which the variable is evaluated.
The convection velocity field $\Vcal$ at time $t_k$ can be approximated by the
backward or forward difference of the displacements,
\begin{nalign}\label{eq-wf12}
  \Vcal(t_k) \approx \incdelta \ub\tl{k}/\dt, \quad \mbox{or} \quad\Vcal(t_k) \approx \incdelta \ub\tl{k + 1}/\dt.
\end{nalign}
Further we employ the abbreviation for the state increments,
\begin{nalign}\label{eq-wf15}
  \incdelta \ub\tl{k+1} \mapsto \ub,\quad \incdelta p\tl{k+1}\mapsto p, \\
  \incdelta \ub\tl{k} \mapsto \bar\ub,\quad \incdelta p\tl{k}\mapsto \bar p\\ 
\end{nalign}
and we refer to the actual spatial configuration by $\Om = \Om\tl{k} = \Om(t_k)$, see
Fig.~\ref{fig-ulf}.
We shall introduce the following tensors (employing the Kronecker symbol $\Ib = (\delta_{ij})$):
\begin{nalign}\label{eq-wf13}
  \Bb(\vb) &= (\nabla\cdot\vb)\Ib - (\nabla \vb)^T,\\
  \Hb(\vb) &= (\nabla\cdot\vb) \Kb - \Kb(\nabla\vb)^T - (\nabla\vb)\Kb^T
\end{nalign}
and the tangent elastic operator
\begin{nalign}\label{eq-wf14}
  \Aop = \Dop^\eff(t) + \sigmabf^\eff(t) \otimes \Ib - p(t) (\Ib\otimes\Ib - \Iop),
\end{nalign}
where $\Iop = (\delta_{jl}\delta_{ik})$. See \ref{app-weak_formulation} for details.

Using expressions \eq{eq-wf07}-\eq{eq-wf12}, the abbreviation \eq{eq-wf15}
and tensors $\Bb$, $\Hb$, $\Aop$ defined in \eq{eq-wf13}, \eq{eq-wf14},
the incremental deformation-diffusion problem represented by \eq{eq-wf03} can be formulated as follows:
Find $(\ub,p) \in \incdelta \Vb \times \incdelta Q$ ($\incdelta V$ and
$\incdelta Q$ are sets of admissible increments) such that
\begin{nalign}\label{eq-wf16}
  \int_{\Om} \Aop \eb(\ub):\eb(\vb)
  - \int_{\Om} p \left(\Bb(\bar\ub) + \Ib \right):\nabla\vb
  =  - \int_{\Om}\left(\sigmabf\tl{k}: \nabla \vb - \fb\tl{k+1}\cdot \vb\right),
\end{nalign}
for all $\vb \in \Vbzero$ and
\begin{nalign}\label{eq-wf17}
  \int_{\Om} &q \left(\Bb(\bar\ub) + \Ib\right):\nabla\ub 
  + \dt \int_{\Om} \left( \Kb + \Hb(\bar\ub)\right)\nabla p \cdot \nabla q\\
  &= \dt \Jcal\tl{k+1}(q) -\dt \int_{\Om}
  \left (\Kb + \Hb(\bar\ub) +  \incdelta \Kb\right)\nabla p\tl{k} \cdot \nabla q,
\end{nalign}
for all $q \in Q_0$. The detailed derivation of equations \eq{eq-wf16}, \eq{eq-wf17} and expressions
\eq{eq-wf13}, \eq{eq-wf14} is in \ref{app-weak_formulation}.
See also our previous work \cite{Rohan-Lukes-ADES2017}, where the similar model,
with respecting inertia effects, is treated.
%
Tensor $\Hb(\vb)$ is symmetric for any symmetric
$\Kb$, whereas $\Bb(\vb)$ is non-symmetric in general. The effective stress
$\sigmabf^\eff$ and its Truesdell rate $\Dop^\eff$ are associated with a strain
energy function and they are functions of deformation. We assume
that $\Dop^\eff$ is elliptical and symmetric.

\section{Homogenization}\label{sec-homogenization}

Using the two-scale homogenization based on the asymptotic analysis of the
deformation-diffusion problem stated above for $\veps \rightarrow 0$, we derive
a limit model describing the medium macroscopic response. by virtue of the
homogenization, the spatial coordinates can be split into the macroscopic parts
denoted by $x$ and the local (microscopic) parts, denoted by $y$. As the
consequence of the finite deformation and the incremental formulation, an
updated configuration is established in terms of the locally representative
periodic deformed cell.

\subsection{Locally periodic microstructures and incremental formulations}

Since the reference configuration is associated with the deformed structure,
this latter is not periodic in general as the result of the nonuniform
deformation. In order to apply the tools of the ``periodic'' homogenization,
the notion of locally periodic structures must be introduced. In fact, the
results of the asymptotic analysis due to the two-scale convergence or periodic
unfolding methods \cite{Cioranescu-etal-UF-book2018} can be obtained even for
microsctructures which vary with the macroscopic position $x \in \Om$ in the
deformed configuration. In particular, the so-called slowly varying
quasi-periodic structures \cite{Brown-Efendiev-IJG2011,Brown-Efendiev-MMS2013}
present a convenient approximate representation of real media. For large
deforming media, the computational algorithm associated with a
time\footnote{The ``time'' $t$ can represent just scalar parameter associated
with iterations of solving a non-linear problem.} discretization $\{t_k\}_k$ is
based on an incremental formulation which enables to determine configurations
at time $t_k+\dt$ based on information available in past, \ie for $t \leq t_k$.

In this paper we consider an incremental two-scale formulation established
using the homogenization of the incremental problem for the ``locally
periodic'' heterogeneous medium. In contrast with linear problems, see \eg
\cite{van-noorden_muntean_2011,Brown-Efendiev-MMS2013}, the local periodicity
property must be updated during time-stepping, or iterations of solving the
nonlinear problem.
For a given scale $\veps$, by $Y^\veps = \veps Y$ refer to the RVE which
coincides with the representative periodic cell $Y$ defined usually as a
parallelepiped; for simplicity, we consider $Y = ]-\frac{1}{2},\frac{1}{2}[^3$
(for the 3D structures) The periodic heterogeneous material is represented by
the microconfiguration $\Mcal^\veps(x,Y,\HH,\SS)$ describing the continuum
occupying domain $Y^\veps(x)$ located at position $x$, whereby its material
(mechanical) properties are given by $\HH^\veps$ comprising all material
parameters, whereas $\SS^\veps$ determines the state variables; both
$\HH^\veps$ and $\SS^\veps$ are functions defined in $Y^\veps$. Note that both
$\HH^\veps$ and $\SS^\veps$ depend on the scale, in general, however, below we
drop the superscript $^\veps$ to lighten the notation.

\paragraph{Local periodicity and the computational procedure}

The assumption of the local periodicity relying on the slowly varying
microconfigurations is needed to homogenize the incremental problem defined by
the linearization and time discretization of the evolutionary problem.
According to the spatial discretization based on the finite element method, the
microconfigurations $\hat \Mcal(\hat x,\hat Y,\hat\HH,\hat\SS)$ are established
at any selected node $\hat x$ associated with the spatial discretization of the
macroscopic problem, such as the integration points $\hat x$, or finite
elements $E(\hat x)$ ``centered'' at $\hat x$. Intuitively, the local
periodicity can be considered within $E(\hat x)$, or in a neighbourhood
$\Ucal_\delta^\veps(\tilde x)$ of the (macroscopic) integration point.

The incremental two-scale formulation for the homogenized medium is related to
the original hetrogeneous one by virtue of some notions explained in the text
below, which are involved in the following time-stepping algorithm:
\begin{algorithm}
  \caption{Time-stepping algorithm}
  \begin{enumerate}[(i)]
      \item Time level $t_k$: Given a consistent heterogeneous medium for a
      given scale $\veps = \veps_0 > 0$, establish an approximate locally
      periodic structure $\Lcal^\veps$, see \eq{eq-3}, in a
      $\delta$-neighbourhood $\Ucal_\delta^\veps(\tilde x)$ of any selected
      node $\tilde x$ associated with the spatial discretization of the
      macroscopic problem.
      \item Using the asymptotic analysis of the incremental problem and the
      local unfolding within the neighbourhood $\Ucal_\delta^\veps(\tilde x)$,
      establish the homogenized two-scale models for all considered $\tilde x$.
      \item For all limit microconfigurations $\Mcal^0(\tilde x,\hat Y,\hat
      \HH,\hat \SS)$, solve local problems for the characteristic responses and compute the
      homogenized coefficients of the macroscopic model.
      \item Collecting the information provided by the microconfigurations $\Mcal^0(\tilde x,\hat Y,\hat
      \HH,\hat \SS)$,
      at all $\tilde x$ considered, constitute the macroscopic numerical model
      of the discretized macroscopic problem which is then solved for the
      macroscopic increments of the state variables, providing updates
      $\Delta\SS$ for each microconfiguration $\Mcal^0$ at time $t_{k+1}
      =t_k+\dt$. Update the local RVEs at $\tilde x$ and time $t_{k+1}$.
      \item Given RVEs at $\tilde x$ (the selected nodes associated with the
      spatial discretization of the macroscopic problem), for the given finite
      scale $\veps_0$, construct an approximate heterogeneous consistent
      structure using the concept of slowly varying microconfiguration.
      \item Updated macroscopic loads, for computing the new time increment,
      with $k:=k+1$ proceed by step (i).
    \end{enumerate}
\end{algorithm}

Although the computational procedure can avoid the reconstruction step (v), an
interpolation procedure should be applied after each time increment computations
(while solving the local problems for all local microconfigurations and the
macroscopic problem) to respect a final scale $\veps_0 > 0$ of the
heterogeneous structure. Indeed, this enables to keep the subsequent
homogenization at time levels $t_1,t_2,\dots$ independent of the spatial
discretization and, thus, to avoid artifacts generated by the discretized
formulation.

\paragraph{Lattice generating a periodic microstructure of the initial configuration}

Let $ \Xi^\veps$ be a set of points in $\RR^3$ and $\vec{a}^j$, $j = 1,2,3$ in
$\RR^3$ be three non-planar vectors $\vec{a}^1\cdot(\vec{a}^2\times\vec{a}^3)
\not = 0$, such that for any two $\xi,\xi'$, there is a $\boldsymbol{k} \in
\ZZ^3$ which yields the path between the two points $\xi' = \xi + \veps k_j
\vec{a}^j$. The periodic (Bravais) lattice $\Lcal_0^\veps(\Xi^\veps,Y^\veps)$
is generated by $\Xi^\veps$ and $\veps Y \equiv Y^\veps$, where $Y^\veps$ is
the representative periodic cell (RPC) defined by the parallelepiped $Y$ with
its three edges constituted by vectors $\vec{a}^k$. By $\Ical^\veps$ we denote
the index set of lattice centers (nodes), \ie $\xi^j \in \Xi^\veps
\Leftrightarrow j \in \Ical^\veps$. We assume that $\Om_0$ is generated by the
periodic lattice $\Lcal^\veps(\xi^0,\Xi^\veps,Y^\veps) =
\Lcal_0^\veps(\Xi^\veps,Y^\veps)$ defined at node $\xi^0 \in \Xi^\veps$, such
that
\begin{equation}\label{eq-1}
\begin{split}  
  \Om_0 = \bigcup_{\xi \in \Xi^\veps} Y^\veps(\xi) \cup \Gamma^\veps\;, \quad Y^\veps(\xi) \cap Y^\veps(\xi') = \emptyset \mbox{ for } \xi \not = \xi' \in \Xi^\veps\;,\\
  Y^\veps(\xi^j) = \xi^0 + \veps(Y + \sum_l^3 k_l^{j} \vec{a}^l)\;,\quad  k_l^{j} \in \ZZ,\; j \in \Ical^\veps\;,
\end{split}
\end{equation}
where $Y^\veps(\xi^j)$ is defined using $\boldsymbol{k}^j$ as the copy of the
scaled RPC $\veps Y$ situated at position $\xi$. Above $\Gamma^\veps$ is the
union of all interfaces $\Gamma_{ij} = \pd Y^\veps(\xi^i) \cap \pd
Y^\veps(\xi^j)$, $i,j \in \Ical^\veps$.

\begin{figure}
  \centering
  \includegraphics[width=0.9\linewidth]{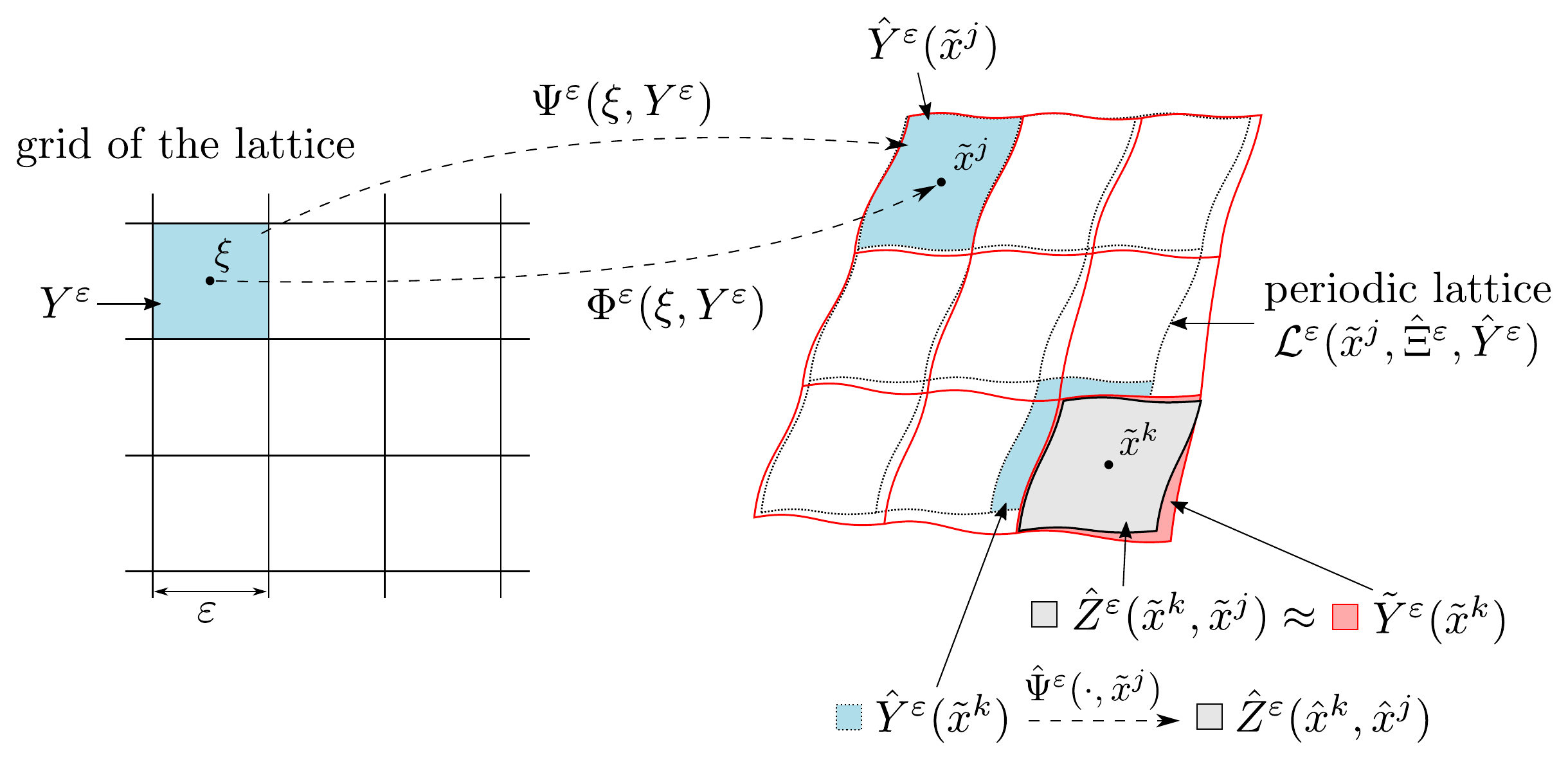}
  \caption{A scheme explaining the  meaning of mappings $\Phi^\veps$ and $\Psi^\veps$ employed to introduce the locally periodic microstructures at the local position $\tilde x^j$ using the lattice $\Lcal^\veps(\hat x^j,\hat\Xi^\veps,\hat Y^\veps)$.}\label{fig-loc-per}
\end{figure}

\paragraph{Slowly varying microstructure and deformed configuration}

The notion of the slowly varying microstructure has been reported \eg in
\cite{Brown-Efendiev-IJG2011} in the context of a stationary Stokes flow
homogenization, where the microstructure is represented by the geometry only.
When dealing with an evolutionary problem of large deforming media, the ``slow
variation'' is influenced by the assumed inhomogeneity of the macroscopic
deformation and of other related state variables which determine the actual
micro-configurations.

We shall assume that the deformed configuration represents a medium with a
slowly varying microstructure, being constituted as the union of cells $\tilde
Y^\veps(x^j) \equiv \tilde Y^{\veps,j}$, such that
\begin{equation}\label{eq-2}
\begin{split}
  \Om & = \bigcup_j \tilde Y^\veps(x^j) \cup \tilde \Gamma^\veps\;,\\
  \tilde Y^\veps(\tilde x) & = \Phi^\veps(\xi,Y^\veps)\;,\quad \tilde x = \vphi^\veps(\xi)\;,\quad \forall \xi \in \Xi^\veps\;,
\end{split}
\end{equation}
where $\Phi^\veps$ is established using a differentiable mapping,
$\vphi^\veps:\Om_0 \mapsto \Om$. As will be explained below, mapping $
\Phi^\veps$ is constructed using the unfolding and a folding procedures
associated with the homogenization $\veps\rightarrow 0$ and the reconstruction
of the homogenized incremental solutions.

\paragraph{Local periodicity and local lattice}

The local periodicity property of the microsctructures is needed to apply the
homogenization procedure; in \cite{Ptashnyk-Loc-per-mi-MMS2015}, a rigorous
treatment was explained in terms of the unfolding approach. We can introduce
locally periodic microstructures as an approximation of slowly varying
microstructures which are assumed to describe the heterogeneous porous medium
characterized by a given finite scale $\veps > 0$, see Fig.~\ref{fig-loc-per}.
In the limit $\veps\rightarrow 0$, due to the perfect scale separation, the
downscaling procedure applied to the homogenized medium yields locally periodic
microstructures. These also naturally represent the numerical approximation of
the homogenized medium, as explained below.

We shall assume that for any $\tilde x = \vphi(\xi)$, such that $\xi \in
\Xi^\veps$, there exist a neighbourhood $\Ucal_\delta^\veps(\tilde x) \subset
\Om$ in which the spatial (deformed) heterogeneous structure can be well
approximated by a locally periodic one. For any $\tilde x^j =
\vphi^\veps(\xi^j)$, $\xi^j \in \Xi^\veps$, there exists the local RPC $\hat Y$
defining the lattice $\Lcal^\veps(\tilde x^j,\hat\Xi,\hat Y)$, such that
\begin{equation}\label{eq-3}
  \begin{split}
    \hat Y^\veps(\hat x^k) & = \veps (\hat Y + \hat z^{k,j})\;,\quad \mbox{ where }  \hat z^{k,j} = \veps^{-1} (\hat x^k - \hat x^j)\;, \; \hat x^j \equiv \tilde x^j\;, \\
    \tilde Z^\veps(\tilde x^k;\tilde x^j) & = \{x \in \Om|\; x = \hat\psi^\veps(\hat x, \tilde x^j),\; \mbox{ for } \hat x \in \hat Y^\veps(\hat x^k)\}\;,\\
\end{split}
\end{equation}
such that $\hat\psi^\veps(\cdot, \tilde x^j)$ is a diffeomorphism mapping the
positions in the local periodic structure on the spatial positions. Hence,
$\tilde Z^\veps(\tilde x^k;\tilde x^j) \approx \tilde Y^\veps(\tilde x^k) =
\Phi^\veps(\xi^k,Y^\veps)$. Below we introduce RPC $\hat Y^\veps(\hat x^j)$ as
an (in a sense, the best) approximation of $\tilde Y^\veps(\tilde x^j)$, see
problem \eq{eq-11}. Consequently, the local lattice $\Lcal^\veps(\hat
x^j,\hat\Xi^\veps,\hat Y^\veps)$ can be defined. Domain $\hat Y^\veps$ can be
considered as the image of a $Y$-periodic regular mapping of cell $\veps Y$.

It should be noted that cells $\tilde Z^\veps(\tilde x^k,\tilde x^j)$ at a
fixed position $\tilde x^k$ can be introduced using different locally periodic
approximations when choosing different $\tilde x^j$. This constitutes the basis
for the homogenization.

Besides the geometry (\ie the RVE decomposition into various compartments
occupied by different phases), the deformed structure is characterized by its
parameters and state variables (fields). Therefore, the approximation property
must respect the ``microconfiguration'', as introduced above. Let
$\tilde\Mcal^\veps(\tilde x,\tilde Y^\veps,\tilde\HH,\tilde\SS)$ be the spatial
configuration associated with domain $\tilde Y^\veps(\tilde x)$. For $\veps >
0$, we can establish an approximation $\hat\Mcal^\veps(\tilde x,\hat
Y,\hat\HH,\hat\SS) \approx \tilde\Mcal^\veps(\tilde x,\tilde
Y^\veps,\tilde\HH,\tilde\SS)$. By $\hat g$ we denote a function defined in the
local lattice $\Lcal^\veps(\hat x,\hat\Xi,\hat Y)$ with $\hat x = \tilde x$ and
$\tilde g^\veps$ the corresponding function defined in the spatial
configuration $\Om$. In general, $\hat g$ can be decomposed into a periodic
part, $\hat g_\#(y,\hat x)$, $y \in \hat Y$, and a macroscopic part, $\hat
g_G(x,\tilde x)$, $x \in Y^{\veps}(\tilde x)$, thus,
\begin{equation}\label{eq-4}
  \begin{split}
    \Tuf{\hat g^\veps(x)} & := \hat g_\#(y,\tilde x) \hat g_G(x,\tilde x)\;,\quad \tilde x = \hat x\;,
\end{split}
\end{equation}
where $\Tuf{}$ is the unfolding operator \cite{Cioranescu-etal-UF-book2018},
cf.~\cite{Ptashnyk-Loc-per-mi-MMS2015}. We assume the following approximation
property is satisfied by $\tilde g^\veps$,
\begin{equation}\label{eq-5}
  \begin{split}
   \|\Tuf{\tilde g^\veps(\hat\psi^\veps(\veps\cdot, \tilde x^j))}- \Tuf{\hat g^\veps}\|_{L^2(\hat Y)} & \rightarrow 0 \quad \mbox{ for } \veps\rightarrow 0\;,
  \end{split}
\end{equation}
where $\veps \cdot$ is substituted by $\veps y$, $y \in \hat Y$.
  
In the generic sense, function $\tilde g^\veps$ represents any material, or
state function of $\tilde\Mcal^\veps(\tilde x,\tilde Y^\veps,\cdot,\cdot)$,
whereas $\hat g^\veps$ is the corresponding approximation.

\paragraph{Reconstruction of continuously varying microstructure for $\veps_0 > 0$}

It is well known, that, in the numerical homogenization, the microstructures are
associated with the discretization scheme. As will be shown below, an
incremental algorithm can be designed which uses the limit two-scale model and
the computation is based on alternating micro- and macro steps. In the context
of such an incremental formulation providing responses at time $t+ \dt$ based
on a known configuration at time $t$, it is desirable to establish a link
between the deformed heterogeneous medium and the homogenized medium for which
the assumption of a locally periodic microstructure is needed. For this there
are two reasons at least: 1) the homogenization of a heterogeneous medium at
time $t$ independently of the numerical approximation, 2) reconstructions of
the homogenized responses for a given scale.

Therefore, we suggest a ``downscaling'' procedure, see step (v) of the
\emph{Algorithm A1}, which enables to reconstruct the heterogeneous structures
based on locally periodic microconfigurations established by virtue of the
homogenization, step (ii) of the \emph{Algorithm A1}. For this, we first need
to introduce the representative periodic cell of the locally periodic
microstructure which represents the actual ``slowly varying'' microstructure
$\tilde Y^\veps(\tilde x^k)$. There exist mappings $\tilde\vphi^\veps$ and
$\Psi^\veps$, such that
\begin{equation}\label{eq-10}
  \begin{split}
    \tilde\vphi^\veps: \tilde Y^\veps(\tilde x^k) & \mapsto \hat Y^\veps(\hat x^k)\;,\\
    \hat Y^\veps(\hat x^k) & = \Psi^\veps(\xi,Y^\veps)\;,\quad \hat x^k = \vphi^\veps(\xi)\;,\quad \xi \in \Xi^\veps\;,
\end{split}
\end{equation}
thus, $\hat Y^\veps$ is the deformed cell $Y^\veps$ at the lattice point $\xi$ of the inital configuration. Mapping $\tilde\vphi^\veps$ is introduced by virtue of the following minimization problem: Find $\Psi^\veps(\xi,\cdot) \in \Hpdb(Y^\veps)$ such that
\begin{equation}\label{eq-11}
  \begin{split}
    \|\Phi^\veps(\xi,\cdot) - \Psi^\veps(\xi,\cdot)\|_{\Lb^2(Y^\veps)}\rightarrow \min.
  \end{split}
\end{equation}
Alternatively, the $\Lb^2(\Ucal_\delta^\veps(\tilde x))$ norm can be employed
together with the lattice $\Lcal^\veps(\tilde x,\hat\Xi^\veps,\hat Y^\veps)$
established by virtue of \eq{eq-3},
\begin{equation}\label{eq-12}
  \begin{split}
    \|\Tuf{\Phi^\veps(\xi,\cdot)} - \Psi(\xi,\cdot,\cdot)\|_{\Ucal_\delta^\veps(\tilde x)\times \hat Y}\rightarrow \min,
  \end{split}
\end{equation}
where scale is fixed $\veps:=\veps_0 > 0$ and $\Psi(\xi,x,y) =
\Tuf{\Psi^\veps(\xi,x)}$ is the two-scale $\hat Y$-periodic function.

For a given $\veps_0$, the characteristic microstructure size $\ell^{\veps_0}
<< h$, where $h>0$ is the perimeter of the finite element, a locally periodic
microstructures, or the slowly varying microstructures can be established using
problem \eq{eq-11}.
\begin{list}{}{}
  \item A) The locally periodic microstructure defined within a subdomain
  $\Om_h(\tilde x^j) \subset \Ucal_\delta^{\veps_0}(\tilde x^j)$ which can be
  associated with a finite element $E_h(\tilde x^j)$ ``located'' at $\tilde
  x^j$ whose the perimeter is $h$.
  \item B) The slowly varying microstructure introduced via an approximation
  $\tilde \Psi^{\veps_0}(\{\tilde x^k\}_{k \in \Ical_h^{\veps_0}},\cdot)$ which
  is constructed as an interpolation of the locally periodic microstructures
  defined at selected points $\{\tilde x^k\}_{k \in \Ical_h^{\veps_0}}$.
\end{list}
Recall that in both the cases, the geometrical transformations related to cells
$\hat Y^\veps$ are employed as the basis for the approximation of
microconfigurations $\tilde\Mcal^\veps(\tilde x,\tilde Y^\veps,\cdot,\cdot)$.

\subsection{Reference cell decomposition and double porosity}

By $Y$ we now abbreviate the local spatial reference cell $\veps^{-1}
Y^{\veps}$, without writing explicitly the macroscopic location. Since our
interests are in porous media with large contrasts in the permeability
coefficients, we shall introduce the following decomposition of $Y$ into the
sectors of primary and dual porosities.

\begin{figure}
  \centering
  \includegraphics[width=0.7\linewidth]{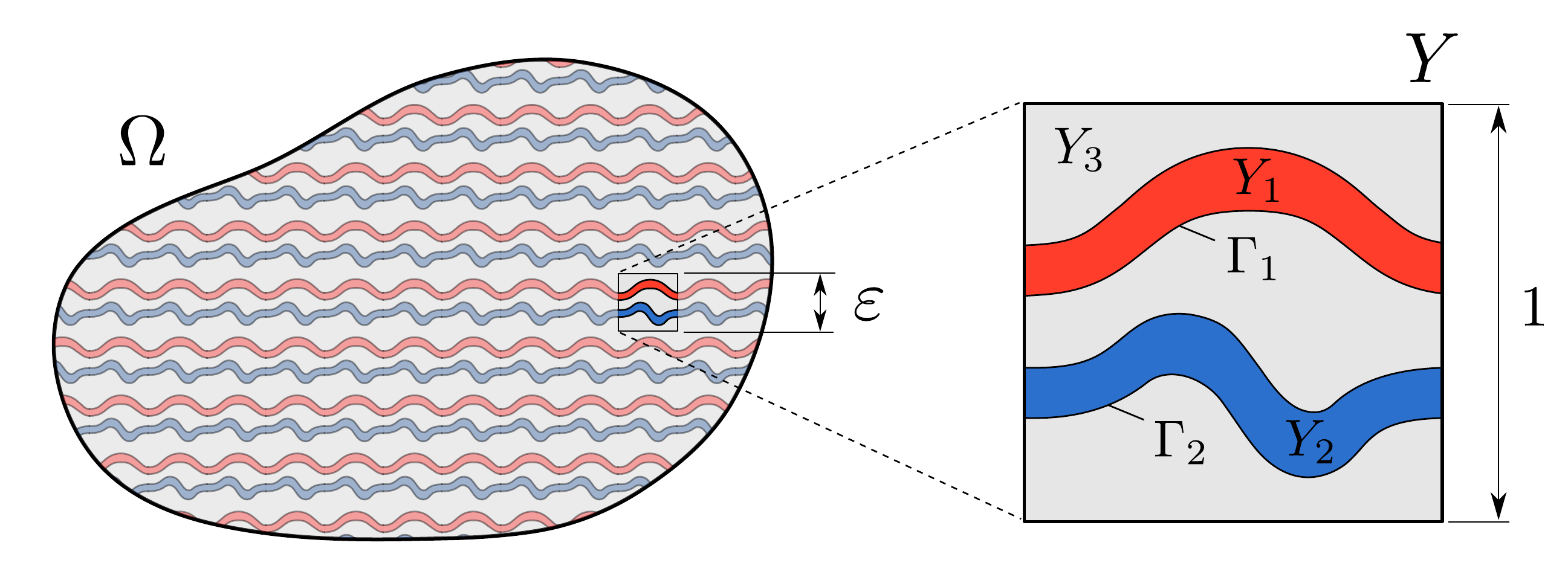}
  \caption{The macroscopic domain $\Om$ and the representative cell
     $Y$ constituted by the porous matrix $Y_3$ and two fluid channels $Y_1$, $Y_2$.}\label{fig-1}
\end{figure}

Let $Y_\alpha$, for $\alpha = 1,2$ be mutually disconnected subdomains of $Y$ with
Lipschitz boundary, then $Y_3$ forms the complement (see Fig.~\ref{fig-1})
\begin{nalign}\label{eq-hom-02}
  Y_3 & \equiv Y \setminus \bigcup_{\alpha = 1,2}\ol{Y_\alpha},\\
  \Gamma_\alpha
    = \pd_3Y_\alpha = \pd_\alpha Y_3 & \equiv \ol{Y_\alpha} \cap \ol{Y_3}.\\
\end{nalign}
Furthermore, we require that $\pd Y \cap \pd Y_l \not = \emptyset$, so that the
all domains $\Om_l^\veps$, $l=1,2,3$ generated by repeating the cells $\veps Y_l$
are connected. The parts of the periodic boundary will be denoted by
\begin{equation}
  \pd_l Y = \pd_l Y_l  \equiv \ol{Y_l} \cap \pd Y\;,\quad l = 1,2,3\;.
\end{equation}

The material parameters in deformed configuration depend on the deformation. To
define them in the local reference cell, we use the unfolding operation,
\cite{Cioranescu2008}. For the sake of brevity, we shall consider couples
$(y,x)$ implicitly matching, such that $y - \hat y = \veps^{-1}(x - \hat x)$.
In the spatial configuration we introduce the permeability parameters, as
follows:
\begin{nalign}\label{eq-hom-03}
  \Kb^\veps(x) = \left \{
  \begin{array}{ll}
  \tilde \Kb^\alpha(x,y) & y \in Y_\alpha, \; \alpha = 1,2, \\
  \veps^2 \tilde \Kb^3(x,y) & y \in Y_3\;.
  \end{array}
  \right .
\end{nalign}
Symbol~$\tilde{}$ denotes quantities that are expressed using the unfolded
spatial coordinates $x$, $y$. We emphasize, that the scaling of the
permeability coefficients by $\veps^2$ in the ``matrix'' compartment leads to
the \emph{double porosity} effect.

\begin{remark}\label{rem-1}
Obviously, in the large deformation regime, tensors $\tilde \Kb^l$ should be
given as functions of the material deformation and, thus, should be modified in
time, as the deformation progresses. By virtue of the homogenization, the local
microscopic deformation gradients are introduced, so that the permeabilities
can be updated at each time level.
\end{remark}

\subsection{Asymptotic expansions and two scale decompositions}

The linearized problem can be treated using the standard homogenization
techniques, such as the periodic unfolding, the two scale convergence, or even
the asymptotic expansion methods.

The displacement and pressure increments $\ub^\veps$ and $p^\veps$ can be
expressed by the following truncated expansions:
\begin{nalign}\label{eq-hom-04}
  \ub^\veps(x) & = \ub^0(x) + \veps \ub^1(x,y) + O(\veps^2),\\
  p^\veps(x) & = \sumalpha\chi_\alpha(y) \left ( p_\alpha^0(x)
    + \veps p_\alpha^1(x,y) + O(\veps^2)\right ) + \chi_3(y) \left (p_3(x,y) + O(\veps) \right ),
\end{nalign}
where we respect the abbreviated notation~\eq{eq-wf15} and $\chi_l$ is the
characteristic function of subdomain $Y_l$ such that $\chi_l(y) = 1$ in $Y_l$
whereas $\chi_l(y) = 0$ in $Y \setminus Y_l$.

Upon substituting these expansions into \eq{eq-wf16}--\eq{eq-wf17} and using
the unfolding technique or the two-scale convergence method, we obtain the
limit problem (for $\veps \rightarrow 0$) formed by the coupled system
involving the unknown displacement increments $(\ub^0, \ub^1) \in \Vbzero(\Om)
\times L^2(\Omega; \HpbY)$ and the pressure increments $(p^0_\alpha,
p^1_\alpha) \in H^1(\Om) \times L^2(\Omega; \Hpalpha))$, $\alpha = 1,2$, $p_3 \in L^2(\Omega;
\Hpzero)$:
\begin{nalign}\label{eq-hom-05}
  \int_{\Om}\intY&\left[\Aop\left(\nabla_x\ub^0 + \nabla_y\ub^1\right)
   - (\sumalpha \chi_\alpha p_\alpha^0 + \chi_3 p_3)(\tilde\Bb(\bar\ub) + \Ib)\right]
   :\left(\nabla_x\vb^0 + \nabla_y\vb^1\right)\\
   =& \int_{\Om}\intY \timek\sigmabf :\left(\nabla_x\vb^0 + \nabla_y\vb^1\right)
    + \int_{\Om}\intY \fb^{new} \cdot \vb^0,
\end{nalign}
\begin{nalign}\label{eq-hom-06}
  \int_{\Om}\intY&(\sumalpha \chi_\alpha q_\alpha^0 + \chi_3 q_3)(\tilde\Bb(\bar\ub) + \Ib)
    :\left(\nabla_x\ub^0 + \nabla_y\ub^1\right)\\
  &+ \dt \int_{\Om}\intYm \left[\tilde\Kb^3 + \tilde\Hb^3(\bar\ub) \right] \nabla_y p_3 \cdot \nabla_y q_3\\
  &+ \dt \int_{\Om} \sumalpha\intYa \left[\tilde\Kb^\alpha + \tilde\Hb^\alpha(\bar\ub) \right]
  (\nabla_x p_\alpha^0  + \nabla_y p_\alpha^1)\cdot (\nabla_x q_\alpha^0 + \nabla_y q_\alpha^1)\\
  =& 
  - \dt \int_{\Om}\intYm \left[\tilde\Kb^3 + \tilde\Hb^3(\bar\ub) + \delta\tilde\Kb^3\right] \nabla_y \timek{p}_3 \cdot \nabla_y q_3\\
  &- \dt \int_{\Om} \sumalpha\intYa \left[\tilde\Kb^\alpha + \tilde\Hb^\alpha(\bar\ub) + \delta\tilde\Kb^\alpha\right]
    (\nabla_x \timek{p}_\alpha^0  + \nabla_y \timek{p}_\alpha^1)\cdot (\nabla_x q_\alpha^0 + \nabla_y q_\alpha^1),
\end{nalign}
where $\nabla_x$ and $\nabla_y$ are the partial gradients with respect to the
macroscopic and microscopic coordinates. The quantities at the right hand sides
used at time $t$, i.e.\ $\sigmabf\tl{k}$, $p_3\tl{k}$, $(p_\alpha^0)\tl{k}$,
$(p_\alpha^1)\tl{k}$, are labelled by symbol $\timek{\ }$ to avoid the double upper
indices in the notation. The average integral for any $Y_l \subset Y$, for the
subdomain index $l = 1,2,3$, is abbreviated by
\begin{nalign}\label{eq-hom-06b}
  \intYl f = \frac{1}{\vert Y\vert}\int_{Y_l} f.
\end{nalign}
Tensors $\tilde\Bb$, $\tilde\Hb^l$ depending on a two-scale vector fields $\vb = (\vb^0, \vb^1)$ are
defined as follows:
\begin{nalign}\label{eq-hom-07}
  \tilde\Bb(\vb) &= (\Ib\otimes\Ib - \Iop)\left(\nabla_x\vb^0 + \nabla_y\vb^1\right),\\
  \tilde\Hb^l(\vb) &= (\nabla_x\cdot\vb^0 + \nabla_y\cdot\vb^1)\tilde\Kb^l
  - \tilde\Kb^l(\nabla_x\vb^0 + \nabla_y\vb^1)^T
  - (\nabla_x\vb^0 + \nabla_y\vb^1)(\tilde\Kb^l)^T.
\end{nalign}
By $H^1(\Om)$ we denote the standard Sobolev space $W^{1,2}(\Om)$, $\Hp(Y)$ is
the functional space of $Y$-periodic functions, $\HpbY$ is the generalization
for vector functions and $\Hpzero$ includes all functions from $\Hp(Y_3)$ that
are zero on the boundary $\Gamma_3$.

Due to the linearity of the limit equations, we can define the decomposition of
the fluctuating functions, introducing the characteristic responses $\omegabf$,
$\pi$, $\eta$ and the particular responses $\ub^P$, $p_l^P$ depending on the
actual pressures and stresses,
\begin{nalign}\label{eq-hom-08}
  \ub^1 & = \omegabf^{ij}\pd_j^x u_i^0 + \sumalpha \omegabf^\alpha p_\alpha^0 + \ub^P,\\
  p_3 & = \pi^{ij}\pd_j^x u_i^0 + \sumalpha \pi^\alpha p_\alpha^0 + p_3^P,\\
  p_\alpha^1 & = \eta_\alpha^i \pd_i^x p_\alpha^0 + p_\alpha^P.
\end{nalign}

\subsection{Local microscopic subproblems}

Further, to simplify the notation, we shall employ the following linear and
bilinear forms:
\begin{nalign}\label{eq-mic-01}
  \alin{\ub}{\vb} & = \intY \Aop \nabla_y\ub :\nabla_y\vb = \intY \Aop \eb_y({\ub}):\eb_y({\vb}),\\
  \blin{l}{q}{\ub} & = \intYl q \left [\tilde\Bb(\bar\ub) + \Ib\right]:\nabla_y\ub,\\
  \clin{l}{p}{q} & = \intYl \left [\tilde \Kb^l + \tilde\Hb^l(\bar\ub)\right]:\nabla_y p \cdot \nabla_y q,\\
  \dlin{l}{p}{q} & = \intYl \delta \tilde\Kb^l \nabla_y p\cdot \nabla_y q,\\
  g_3^\alpha(\bar\ub;p) & = \int_{\Gamma_\alpha} \left( \tilde\Kb^3
    + \tilde\Hb(\bar\ub)\right) \nabla_y p \cdot \nb^3 \dSy,
\end{nalign}
defined for $l = 1,2,3$ and for $\alpha=1,2$. By $\nb^3$ we denote the normal
unit vector on boundary $\Gamma_\alpha = \ol{Y_\alpha} \cap \ol{Y_3}$ oriented
outwards of $Y_3$. Note that all integrals are defined in the reference
deformed configuration. 

The microscopic equations are obtained from \eq{eq-hom-05}--\eq{eq-hom-06} by
letting vanish all the macroscopic test functions,
\begin{nalign}\label{eq-mic-02}
  \alin{\Pibf:\nabla_x\ub^0 + \ub^1}{\vb}
  - \sumalpha\blin{\alpha}{p^0_\alpha}{\vb}
  - \blin{3}{p_3}{\vb}
  = - \intY \timek\sigmabf:\nabla_y \vb,\\
  \blin{3}{q}{\Pibf:\nabla_x\ub^0 + \ub^1} + \dt \clin{3}{p_3}{q}
  = -\dt\left[\clin{3}{\timek{p}_3}{q} + \dlin{3}{\timek{p}_3}{q}\right],
\end{nalign}
holding for all $\vb \in \HpbY$ and $q \in \Hp(Y_3)$, and
\begin{nalign}\label{eq-mic-02b}
  \clin{\alpha}{\yb \cdot \nabla_x p_\alpha^0 + p_\alpha^1}{q}
  =
  -\clin{\alpha}{\yb \cdot \nabla_x \timek{p}_\alpha^0 + \timek{p}_\alpha^1}{q}
  -\dlin{\alpha}{\yb \cdot \nabla_x \timek{p}_\alpha^0 + \timek{p}_\alpha^1}{q},
\end{nalign}
for all $q \in \Hpalpha$. Above, we use tensor $\Pibf = (\Pi^{rs}_i)$ with
components $\Pi^{rs}_i = y_s \delta_{ir}$ that involve the microscopic
coordinates $\yb$.

Employing the two-scale decompositions \eq{eq-hom-08} in equations
\eq{eq-mic-02}--\eq{eq-mic-02b} enables us to extract the local problems for the
characteristic responses $\omegabf,\pi$ and $\eta$. The following
problems associated with the poroelasticity in the whole cell $Y$ and perfusion
in the dual porosity $Y_3$ must be resolved for a given time step increment
$\dt$:
\begin{enumerate}
\item
Find $(\omegabf^{ij},\pi^{ij}) \in \HpbY \times \Hpzero$ such that
\begin{nalign}\label{eq-mic-03}
  \alin{\omegabf^{ij} + \Pibf^{ij}}{\vb} - \blin{3}{\pi^{ij}}{\vb} & = 0,
    \quad &\forall \vb \in \HpbY,\\
  \blin{3}{q}{\omegabf^{ij} + \Pibf^{ij}} + \dt \clin{3}{\pi^{ij}}{q} & = 0,
    &\forall q \in \Hpzero.
\end{nalign}
\item
Find $(\omegabf^\alpha,\pi^\alpha) \in \HpbY \times  H_{\#}^1(Y_3)$  such that
\begin{nalign}\label{eq-mic-04}
  \alin{\omegabf^\alpha}{\vb} - \blin{3}{\pi^\alpha}{\vb} & = \blin{\alpha}{1}{\vb},
    \quad &\forall \vb \in \HpbY, \\
  \blin{3}{q}{\omegabf^\alpha} + \dt \clin{3}{\pi^\alpha}{q} & = 0,
    &\forall q \in \Hpzero,
\end{nalign}
where (in the sense of traces) $\pi^\alpha = \delta_{\alpha\beta}$ on $\Gamma_\beta$.
\item
Find $(\ub^P,p_3^P) \in \HpbY \times \Hpzero$, the particular response to the current reference state,  such that
\begin{nalign}\label{eq-mic-05}
  \alin{\ub^P}{\vb} - \blin{3}{p_3^P}{\vb} & = - \intY\timek\sigmabf : \nabla_y \vb,
    \quad &\forall \vb \in \HpbY, \\
  \blin{3}{q}{\ub^P} + \dt \clin{3}{p_3^P}{q} & = 
  -\dt\left(\clin{3}{\timek{p}_3}{q} + \dlin{3}{\timek{p}_3}{q}\right),
     &\forall q \in \Hpzero.
\end{nalign}

\end{enumerate}
The following two microscopic subproblems are related to the fluid channels
$Y_\alpha$, $\alpha =1,2$:
\begin{enumerate}
  \item The channel flow correctors:
  Find $\eta_\alpha^i \in \Hp(Y_\alpha)$ such that
  \begin{nalign}\label{eq-mic-06}
    \clin{\alpha}{\eta_\alpha^i + y_i}{q} = 0, \quad
      \forall q \in \Hp(Y_\alpha). \\
  \end{nalign}
  
  \item The particular response for the current load response:
  Find $p_\alpha^P \in \Hp(Y_\alpha)$, such that
  \begin{nalign}\label{eq-mic-07}
    \clin{\alpha}{p_\alpha^P}{q}
      = \dlin{\alpha}{\yb\cdot \nabla_x \timek p_\alpha^0 + \timek p_\alpha^1}{q}
      -\clin{\alpha}{\yb\cdot \nabla_x \timek p_\alpha^0 + \timek p_\alpha^1}{q},\quad
    \forall q \in \Hp(Y_\alpha).
  \end{nalign}
\end{enumerate}

\subsection{Global macroscopic problem and the homogenized coefficients}\label{sec-macro}

The global macroscopic equations and the homogenized coefficients are obtained
from system \eq{eq-hom-05}--\eq{eq-hom-06} applying the macroscopic test
functions $\vb^0$, $q^0_\alpha$ and taking $\vb^1 = 0$, $q^1_\alpha = 0$, but
$q_3 = 0$ in $Y_3$ and $q_3 = q^0_\alpha$ on $\Gamma_\alpha$. The balance of
forces equation attains the following form:
\begin{nalign}\label{eq-mac-01}
    \int_{\Om}\intY&\left[\Aop\nabla_y(\Pibf:\nabla_x \ub^0 + \ub^1)
    - \left(\sumalpha \chi_\alpha p_\alpha^0 + \chi_3 p_3\right)(\tilde\Bb(\bar\ub) + \Ib)\right]
    :\nabla_y(\Pibf:\nabla_x \vb^0)\\
    =& \int_{\Om}\intY \timek\sigmabf : \nabla_x\vb^0
    + \int_{\Om}\intY \fb^\new \cdot \vb^0
\end{nalign}
and the diffusion equation results in two limit equations for $\alpha = 1, 2$
\begin{nalign}\label{eq-mac-02}
  \int_{\Om}\intYa& q_\alpha^0 \left[\tilde\Bb(\bar\ub) + \Ib\right]
    :\nabla_y(\Pibf:\nabla_x\ub^0 + \ub^1)
   + \dt \int_{\Om} q^0_\alpha g_3^\alpha(\bar\ub; p_3)\\
  &+ \dt \int_{\Om} \intYa \left[\tilde\Kb^\alpha + \tilde\Hb^\alpha(\bar\ub) \right]
   \nabla_y(\yb\cdot\nabla_x p_\alpha^0 + p_\alpha^1) \cdot \nabla_x q_\alpha^0\\
  =& - \dt \int_{\Om} q^0_\alpha \left(g_3^\alpha(\bar\ub; \timek p_3)
    + \int_{\Gamma_\alpha} \delta\tilde\Kb^3 \nabla_y \timek p_3 \cdot \nb^3 \dSy\right)\\
  &- \dt \int_{\Om} \intYa \left[\tilde\Kb^\alpha + \tilde\Hb^\alpha(\bar\ub) + \delta\tilde\Kb^\alpha\right]
    \nabla_y(\yb \cdot \nabla_x \timek{p}_\alpha^0 + \timek{p}_\alpha^1)\cdot \nabla_y(\yb \cdot \nabla_x q_\alpha^0),
\end{nalign}
where we use the linear form
\begin{nalign}\label{eq-mac-02b}
  g_3^\alpha(\bar\ub; p) = \int_{\Gamma_\alpha} \left(\tilde\Kb^3 + \tilde\Hb(\bar\ub)\right)
    \nabla_y p \cdot \nb^3 \dSy.
\end{nalign}

Upon substituting the microscopic fluctuations $\ub^1$ and $p_3$ in the
equilibrium equation \eq{eq-mac-01} by the decompositions defined in
\eq{eq-hom-08}, the following homogenized coefficients can be introduced:
\begin{itemize}
  \item The effective viscoelastic incremental tensor,  $\Dcalbf = (\Dcal_{ijkl})$,
  \begin{nalign}\label{eq-mac-03}
    \Dcal_{ijkl} & = |Y|^{-1} \left[\alin{\omegabf^{kl} + \Pibf^{kl}}{\Pibf^{ij}}
      - \blin{3}{\pi^{kl}}{\Pibf^{ij}}\right] \\
    & = |Y|^{-1} \left[\alin{\omegabf^{kl} + \Pibf^{kl}}{\omegabf^{ij} + \Pibf^{ij}}
      + \dt \clin{3}{\pi^{kl}}{\pi^{ij}}\right].
  \end{nalign}
  The above symmetric expression is derived in \ref{app-symcf}.

  \item The Biot poroelasticity tensor, $\Bcalbf = (\Bcal_{ij})$,
  \begin{nalign}\label{eq-mac-04}
    \Bcal_{ij}^\alpha & = |Y|^{-1} \left[\blin{3}{\pi^\alpha}{\Pibf^{ij}}
      + \blin{\alpha}{1}{\Pibf^{ij}} - \alin{\omegabf^\alpha}{\Pibf^{ij}}\right].
  \end{nalign}
  
  \item The averaged Cauchy stress, $\Scalbf = (\Scal_{ij})$,
  \begin{nalign}\label{eq-mac-05}
    \Scalbf & = |Y|^{-1} \int_Y \timek{\sigmabf}.
  \end{nalign}
  
  \item The retardation stress, $\Qcalbf = (\Qcal_{ij})$,
  \begin{nalign}\label{eq-mac-06}
    \Qcal_{ij} & = |Y|^{-1} \left[\alin{\ub^P}{\Pibf^{ij}} - \blin{3}{p_3^P}{\Pibf^{ij}}\right].
  \end{nalign}
\end{itemize}
From the diffusion equation \eq{eq-mac-02}, upon substituting the microscopic
functions $p^1_\alpha$ and $p_3$, we can identify the following homogenized
coefficients:
\begin{itemize}
  \item The effective channel permeability, $\Ccalbf = (\Ccal_{ij})$,
  \begin{nalign}\label{eq-mac-07}
    \Ccal_{ij}^\alpha & = |Y|^{-1} \clin{\beta}{\eta^i_\alpha + y_i}{\eta^j_\alpha + y_j}.
  \end{nalign}

  \item The adjoint Biot poroelasticity tensor, $\Rcalbf = (\Rcal_{ij})$,
  \begin{nalign}\label{eq-mac-08}
    \Rcal_{ij}^\alpha = |Y|^{-1} \left[\blin{\alpha}{1}{\omegabf^{ij}+\Pibf^{ij}}
      + \dt g_3^\alpha(\bar\ub,\pi^{ij})\right].
  \end{nalign}

  \item The perfusion coefficient -- inter-channel permeability,
  \begin{nalign}\label{eq-mac-09}
    \Gcal_\beta^\alpha = |Y|^{-1} \left[g_3^\alpha(\bar\ub,\pi^\beta) +
      (\dt)^{-1}\blin{\alpha}{1}{\omegabf^\beta}\right].
  \end{nalign}

 \item The effective discharges due to deformation of the reference state:
  \begin{nalign}\label{eq-mac-10}
    \zeta_\alpha^\effx & = |Y|^{-1}[\dt^{-1}\blin{\alpha}{1}{\ub^P}
      + g_3^\alpha(\bar\ub,p_3^P + \timek p_3) +
      \int_{Y_\alpha} \delta \tilde \Kb^3  \nabla_y \timek p_3 \cdot  \nb^3 \dSy],\\
    \gamma_\alpha^{i \effx}& =  |Y|^{-1}\left[\clin{\alpha}{\yb\cdot \nabla_x \timek{p}_\alpha^0
      + \timek p_\alpha^1}{y_i} + \dlin{3}{\yb\cdot \nabla_x \timek{p}_\alpha^0 + \timek{p}_\alpha^1}{y_i}{\alpha}
      + \clin{\alpha}{p_\alpha^P}{y_i}\right].
  \end{nalign}
\end{itemize}
In the above expressions, the boundary integrals on the interface
$\Gamma_\alpha$, $\alpha=1,2$, can be computed using the residual form
expression, which yields
\begin{nalign}\label{eq-mac-11}
  g_3^\alpha(\bar\ub,\pi^{ij}) & =
    (\dt)^{-1}\blin{3}{\pi^\alpha}{\omegabf^{ij} + \Pibf^{ij}}
    + \clin{3}{\pi^{ij}}{\pi^\alpha}, \\
  g_3^\alpha(\bar\ub,\pi^\beta) & =
    (\dt)^{-1}\blin{3}{\pi^\alpha}{\omegabf^\beta}
    + \clin{3}{\pi^\alpha}{\pi^\beta}, \\
  g_3^\alpha(\bar\ub,p_3^P+ \timek p_3) & =
    (\dt)^{-1}\blin{3}{\pi^\alpha}{\ub^P}
    + \clin{3}{p_3^P + \timek p_3}{\pi^\alpha} + \dlin{3}{\timek p_3}{\pi^\alpha}.
\end{nalign}

The coefficients $\Bcalbf^\alpha$ and $\Rcalbf^\alpha$ couple the macroscopic displacement
field and the pressure fields in two fluid channels. In \ref{app-symcf} we show that
$\Bcalbf^\alpha = \Rcalbf^\alpha$ which leads to the symmetry of the resulting
macroscopic system.

We substitute the homogenized coefficients into the limit equations
\eq{eq-mac-01}--\eq{eq-mac-02}; hence, we obtain the global macroscopic
problem. The macroscopic solution $\ub^0 \in \Vb(\Om)$ and
$p^0_\alpha \in {Q^\alpha}(\Om)$, $\alpha=1,2$ must satisfy the equilibrium
equation
\begin{nalign}\label{eq-mac-12}
  \int_\Om \left(\Dcalbf \nabla_x \ub^0
    - \sumalpha p^0_\alpha \Bcalbf^\alpha\right):\nabla_x \vb
    =  L^{new}(\vb) - \int_\Om \left(\Scalbf + \Qcalbf\right):\nabla_x \vb,
\end{nalign}
for all $\vb \in \Vbzero(\Om)$ and the diffusion equations
\begin{nalign}\label{eq-mac-13}
  \int_\Om q_\alpha \left(\Bcalbf^\alpha :\nabla_x \ub^0
    + \dt \sum_{\beta=1,2} \Gcal_\beta^\alpha p^0_\beta\right)
    + \dt \int_\Om \Ccalbf^\alpha \nabla_x p^0_\alpha \cdot \nabla_x q_\alpha
    = - \dt \int_\Om \left(\zeta_\alpha^\effx q_\alpha
    + \gammabf_\alpha^\effx\cdot \nabla_x q_\alpha \right),
\end{nalign}
for all $q_\alpha \in {Q_0^\alpha}(\Om)$ and for a given $L^\new(\vb)$ which
involves the volume and surface traction forces at time $t + \dt$.
Keep in mind that $\ub^0$ and $p^0_\alpha$ are the increments,
within the context of the Eulerian formulation, of the macroscopic fields.
The admissibility sets $Q^\alpha$ and spaces $Q_0^\alpha$ reflect the boundary
conditions for the pressure increments on boundary $\pd\Om$. In
particular, pressure $p_\alpha$ can be prescribed on $\pd_{p,\alpha}\Om \subset
\pd \Om$, whereas there is no relationship between $\pd_{p,1}\Om$ and
$\pd_{p,2}\Om$, they can be disjoint or can overlap.

\subsection{Updating incremental algorithm}\label{sec-incremental-algorithm}

We shall now explain the time stepping algorithm which is used to compute
deformation of the heterogeneous structure and fluid redistribution in the
pores at discrete time levels. The procedures involved in the coupled
micro-macro computation process are described in Algorithm~\ref{alg1} where we
use again the unabbreviated notation (abbreviation introduced in \eq{eq-wf15},
so that symbol $\incdelta$ denotes the time increments and the upper index
$(k)$ refers to a certain time level. So, in particular,
$(\incdelta\ub^0)\tl{k+1}$, $(\incdelta p^0_{\alpha})\tl{k+1}$ correspond to
$\ub^0$, $p^0_{\alpha}$ in \eq{eq-mac-12}, \eq{eq-mac-13}, see abbreviation
\eq{eq-wf15}. Also note that we have different local configurations
$Y\tl{k}(x)$ for different macroscopic points $x\in\Om\tl{k}$ at a given time
level. In practice, we solve the local subproblems only in a finite number of
the macroscopic points which are usually associated with the quadrature points
of the finite element discretization at the macroscopic scale. Alternatively we
can assume the constant values of the effective parameters within the
macroscopic finite element and reduce the total number of local problems to the
number of mesh elements at the global level.

\begin{algorithm}
\caption{Calculate new microscopic and macroscopic configurations at time level $k + 1$}\label{alg1} 
\begin{algorithmic}
  \REQUIRE Microscopic and macroscopic configurations at time level $k$
  \WHILE{$k \leq k_{max}$}
    \STATE
     \begin{itemize}
        \item For a given reference two-scale configuration solve the microscopic
        subproblems \eq{eq-mic-03}--\eq{eq-mic-07} and evaluate homogenized
        coefficients $\Dcalbf$, $\Bcalbf^\alpha$, $\Scalbf$, $\Qcalbf$, $\Ccalbf^\alpha$,
        $\Gcal^\alpha_\beta$, $\zeta_\alpha$, $\gammabf_\alpha$
        using expressions \eq{eq-mac-02}--\eq{eq-mac-11}.
        \item Compute new macroscopic increments $(\incdelta\ub^0)\tl{k+1}$
          and $(\incdelta p^0_{\alpha})\tl{k+1}$ by solving equations
        \eq{eq-mac-12}, \eq{eq-mac-13} and update macroscopic fields:
          \begin{nalign}\nonumber
            (\ub^0)\tl{k+1} &= (\ub^0)\tl{k} + (\incdelta\ub^0)\tl{k+1},\\
            (p^0_\alpha)\tl{k+1} &= (p^0_\alpha)\tl{k} + (\incdelta p^0_\alpha)\tl{k+1}.
          \end{nalign}
          (Beware that the increments are now written in the unabbreviated notation).
        \item Update macroscopic configuration:
          $\Om\tl{k+1} = \Om\tl{k} + \lbrace (\incdelta\ub^0)\tl{k + 1}\rbrace$.
        \item Update microscopic configurations for selected macroscopic points $x \in \Om\tl{k}$:
          \begin{nalign}\nonumber
            \incdelta\ub & = (\Pibf^{ij} + \omegabf^{ij})\pd_j^x (\incdelta u_i^0(x))\tl{k+1}
               + \sumalpha \omegabf^\alpha (\incdelta p_\alpha^0(x))\tl{k+1}+ \ub^P,\\
            p_3\tl{k+1} & = p_3\tl{k}
              + \pi^{ij}\pd_j^x (\incdelta u_i^0(x))\tl{k+1}
              + \sumalpha \pi^\alpha (\incdelta p_\alpha^0(x))\tl{k+1} + p_3^P,\\
            p_\alpha\tl{k + 1} & = p_\alpha\tl{k} + (\incdelta p_\alpha^0(x))\tl{k+1}
              + \eta_\alpha^i \pd_i^x (\incdelta p_\alpha^0(x))\tl{k+1} + p_\alpha^P.
          \end{nalign}
        \item Update deformed microscopic domains:
          $Y\tl{k+1}(x) = Y\tl{k}(x) + \lbrace \incdelta\ub\rbrace$.
          Consequently, the local microscopic deformation gradients $\Fb_y =
          (F_{ij}(y))$ are computed, $(F_{ij}(y) = \pd{y}{y^0}$, where $y \in
          Y\tl{k+1}(x)$ are the updated spacial coordinates of the material
          point $y^0$.
        \item Recalculate deformation dependent tangential stiffness tensor $\Dop^\eff$ and 
          effective Cauchy stress tensor $\sigmabf^\eff$.
          Optionally, if the permeabilities $\tilde \Kb^l$ are given functions of the micro-level deformation $\Fb_y^T\Fb_y$, the respective values should bue updated, see Remark~\ref{rem-1}.
        \item $k = k + 1$
      \end{itemize}
  \ENDWHILE
\end{algorithmic}
\end{algorithm}


\section{Numerical simulations}\label{sec-num-simul}

The aim of this section is twofold. Firstly, we present a validation test of
the proposed two-scale homogenized model of a double porosity large deforming
poroelastic medium. Secondly, we demonstrate the applicability of the model in
simple illustrative 2D simulations in Sections~\ref{sec:num_shear_test} and \ref{sec:num_infl_test}.

We employ the finite element method for space discretization at both the scales
and we approximate the displacement and pressure fields by piecewise linear
functions. The mathematical model has been implemented in the software called
{\it SfePy} -- Simple Finite Elements in Python, whose development is focused
on multiscale simulations of heterogeneous structures, see
\cite{Cimrman-Lukes-SfePy2019}. The code is written in Python, partially in C,
and it uses external libraries for solving sparse linear systems such as PETSc,
UMFPACK and MUMPS, that ensures high computational performance. The key feature
of {\it SfePy} is the so-called homogenization engine which allows to solve the
local microscopic subproblems and to evaluate the homogenized coefficients in
an efficient way employing either multithreading or multiprocessing
capabilities of a computer system.

\subsection{Validation test}\label{sec:validation-test}

A direct finite element simulation of the non-homogenized periodic structure is
performed in order to calculate a reference solution which is compare with the
results obtained by the two-scale numerical simulation. 

\begin{figure}
  \centering
    \includegraphics[width=0.6\linewidth]{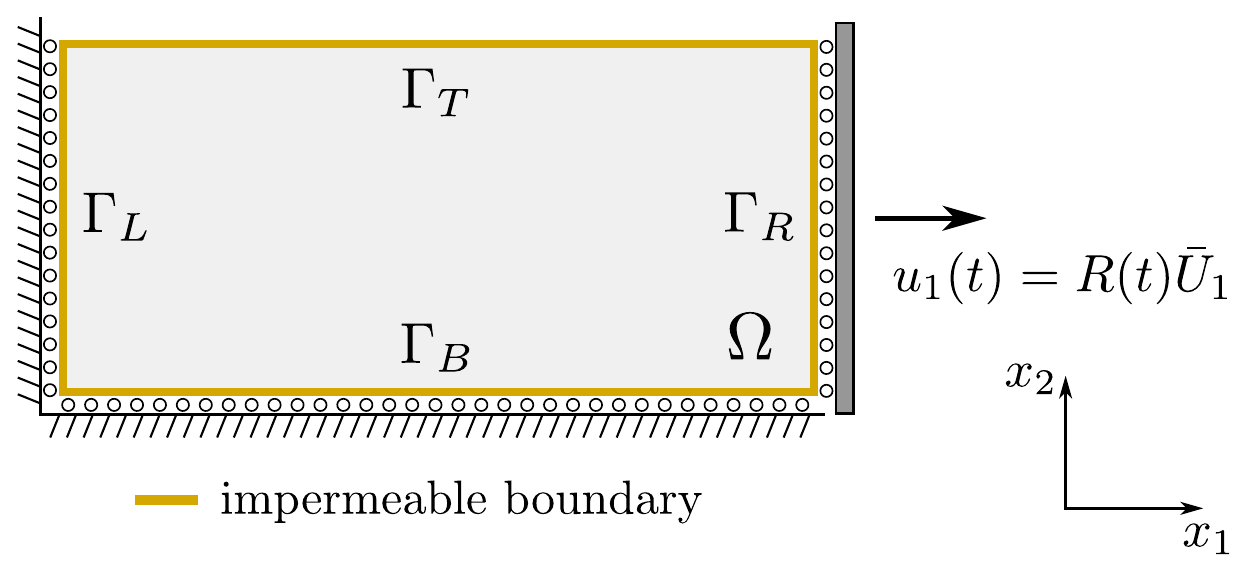}\\
    \includegraphics[width=0.6\linewidth]{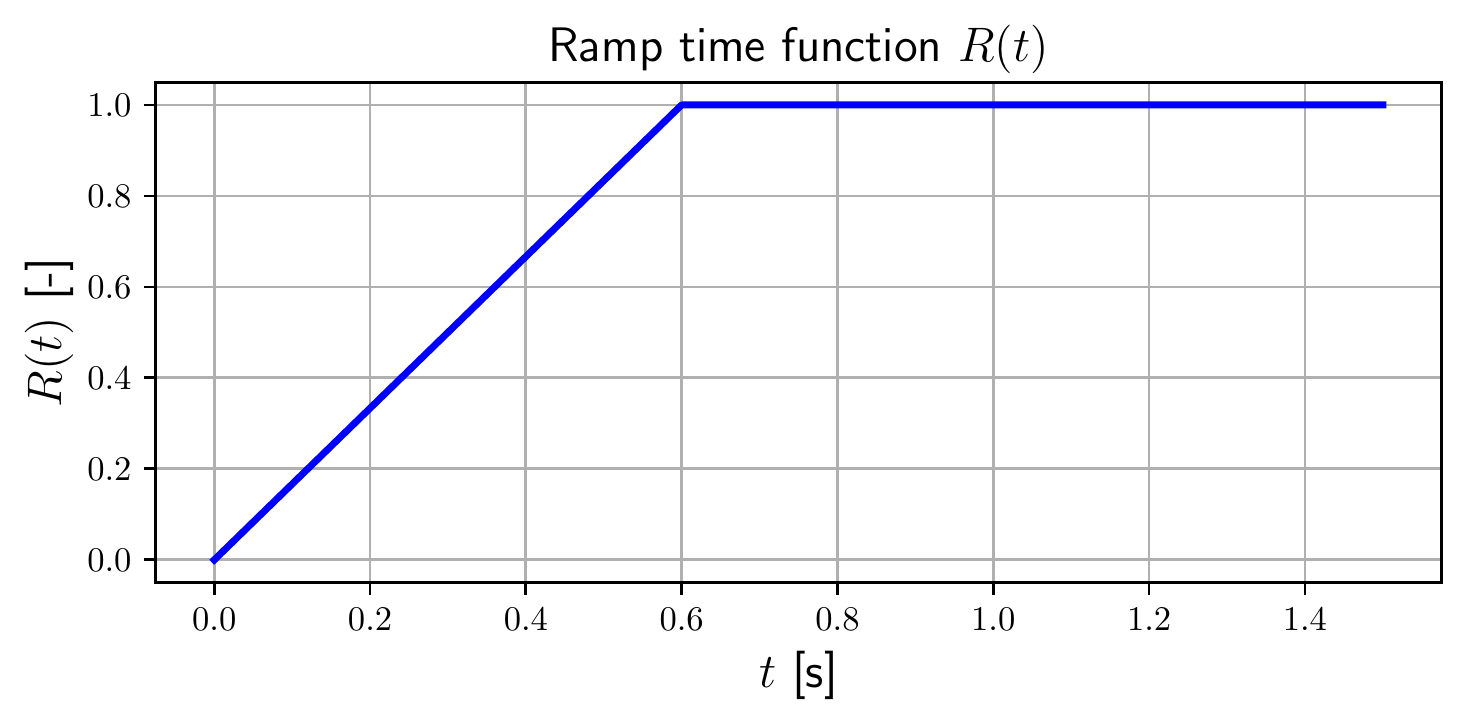}
  \caption{Validation test:
    top -- boundary conditions applied to the macroscopic 2D sample;
    bottom -- ramp function $R(t)$ driving deformations of the macroscopic sample in time.
  }\label{fig-validation_bc_Y}
\end{figure}

\begin{figure}
  \centering
    \includegraphics[width=0.65\linewidth]{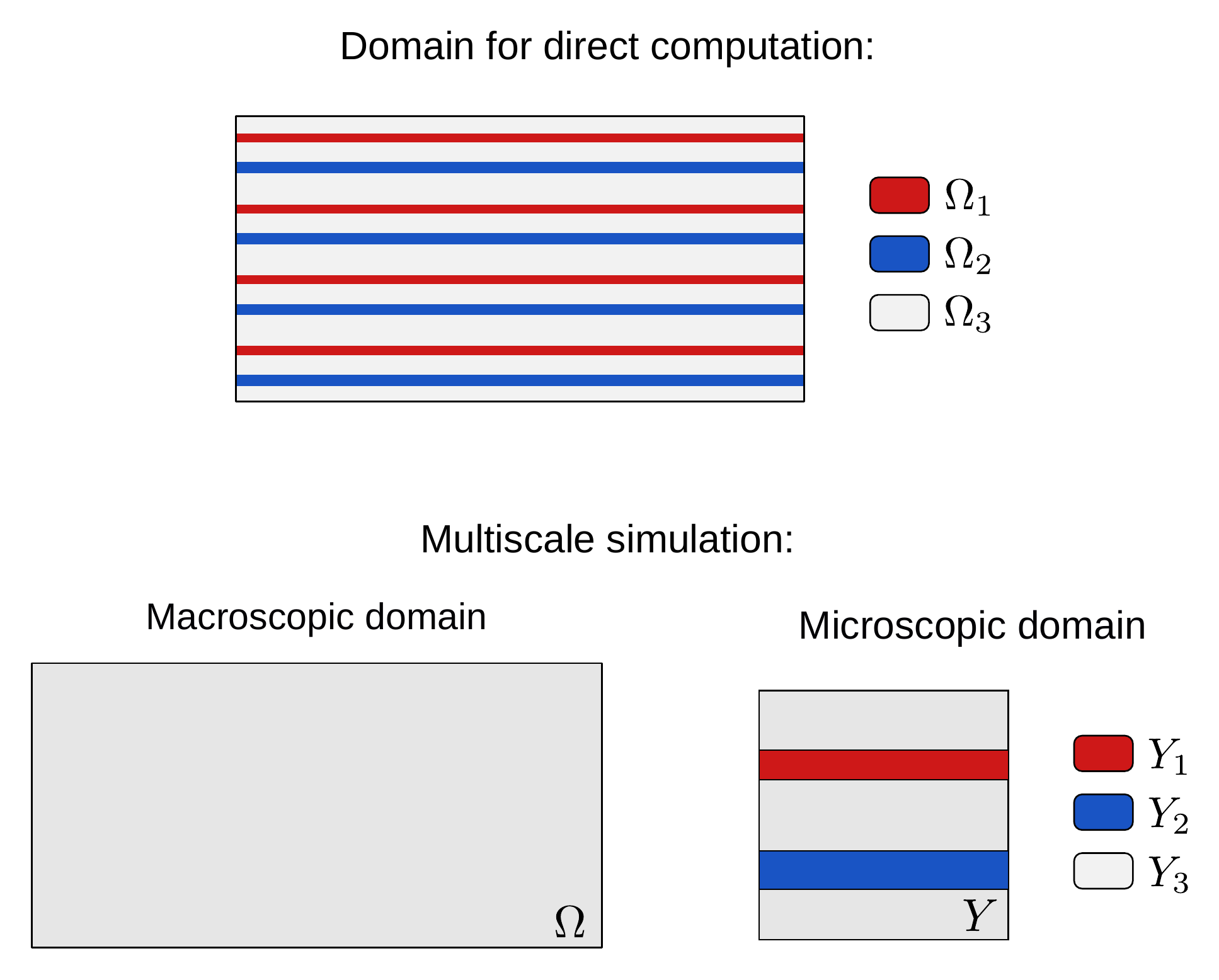}
  \caption{Computational domains:
    top -- domain used in the reference model, build up by repetition ($8\times 4$)
    of the microscopic unit;
    bottom -- macroscopic domain $\Omega$ (left) and decomposition of microscopic 
    domain $Y$ (right).}\label{fig-validation_geometry}
\end{figure}

We consider the 2D rectangular sample with dimensions $0.2 \times 0.1$\,m which
is attached to the rigid frame in such a way that only uniform extension or
contraction in the whole domain is allowed, see Fig.~\ref{fig-validation_bc_Y}
top. The first component of the displacement at the right boundary
($\Gamma_{R}$) is driven in time by function $u_1(t) = R(t) \bar{U}_1$, where
$R(t)$ is the ramp function depicted in Fig.~\ref{fig-validation_bc_Y} bottom and
$\bar{U}_1 = 0.04$\,m. No volume and surface forces are considered. The whole
boundary of the sample is assumed to be impermeable ($\frac{\pd p_\alpha}{\pd
n} =0$ on $\Gamma$). This condition corresponds to closed fluid pores on the outer
surface. At the microscopic level, the representative volume element is
constituted by the hyperelastic solid skeleton and by two straight fluid
channels, see Fig.~\ref{fig-validation_geometry} bottom. The non-linear material
behavior of the solid part is governed by the neo-Hookean constitutive law and
it is determined by the shear modulus $\mu$. The prescribed hydraulic
permeability coefficients $\tilde \Kb^l$ and the elasticity parameters $\mu^l$
in domains $Y_l$, $l=1,2,3$ are summarized in Tab.~\ref{tab-mat_params}. Note
that the permeability in the porous matrix $Y_3$ is rescaled by $\veps^2$, see
\eq{eq-hom-03}. Also note that we need to define some stiffness in the channel
parts, otherwise we get the disconnected structure due to the 2D topology of
the reference cell.

The reference model is built up by copies of the unit periodic cell $Y$. In our
validation test, we take the grid of $8 \times 4$ unit cells rescaled to the
above defined sample dimensions, see Fig.~\ref{fig-validation_geometry} top.
With regard to the sample size and the number of repetitions of the reference
cells, we define the scales ration, employed in the multiscale simulations, as
$\veps=0.2 / 8 = 0.025$.

\begin{table}
  \centering
  \begin{tabular}{|c|c|c|}
    \hline
     & shear modulus & hydraulic permeability\\
     $Y^l$ & $\mu^l$ [MPa] & $\tilde \Kb^l$ [m$^2$ / (Pa s)]\\
    \hline
    channel $Y_1$ & $0.6$ & $10^{-6} \Ib$ \\
    channel $Y_2$ & $0.6$ & $2 \cdot 10^{-6} \Ib$\\
    matrix $Y_3$ & $1$ & $10^{-4} \Ib$ \\
    \hline
  \end{tabular}
  \caption{Permeability and elasticity parameters in
           the porous matrix $Y_3$ and two channels $Y_1$, $Y_2$.}\label{tab-mat_params}
\end{table}

\begin{figure}
  \centering
    \includegraphics[width=0.48\linewidth]{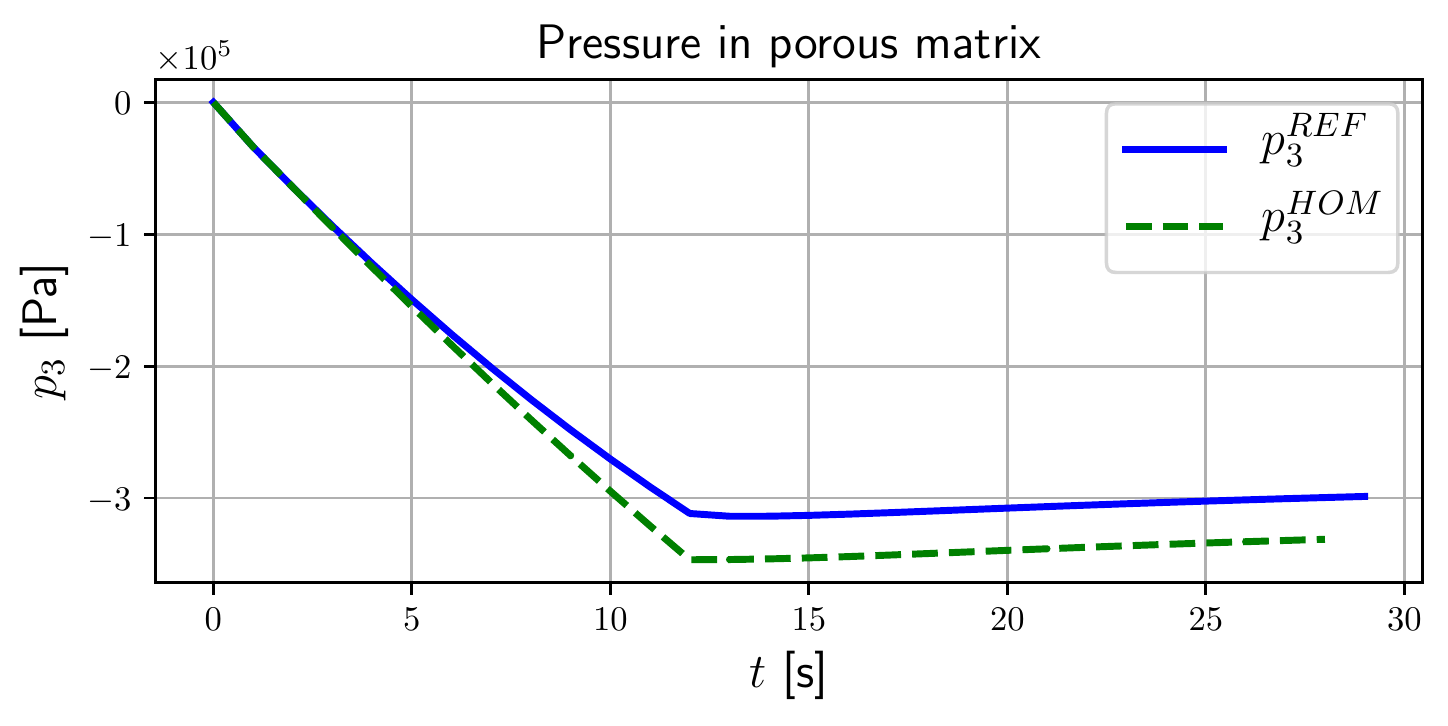}\hfil
    \includegraphics[width=0.48\linewidth]{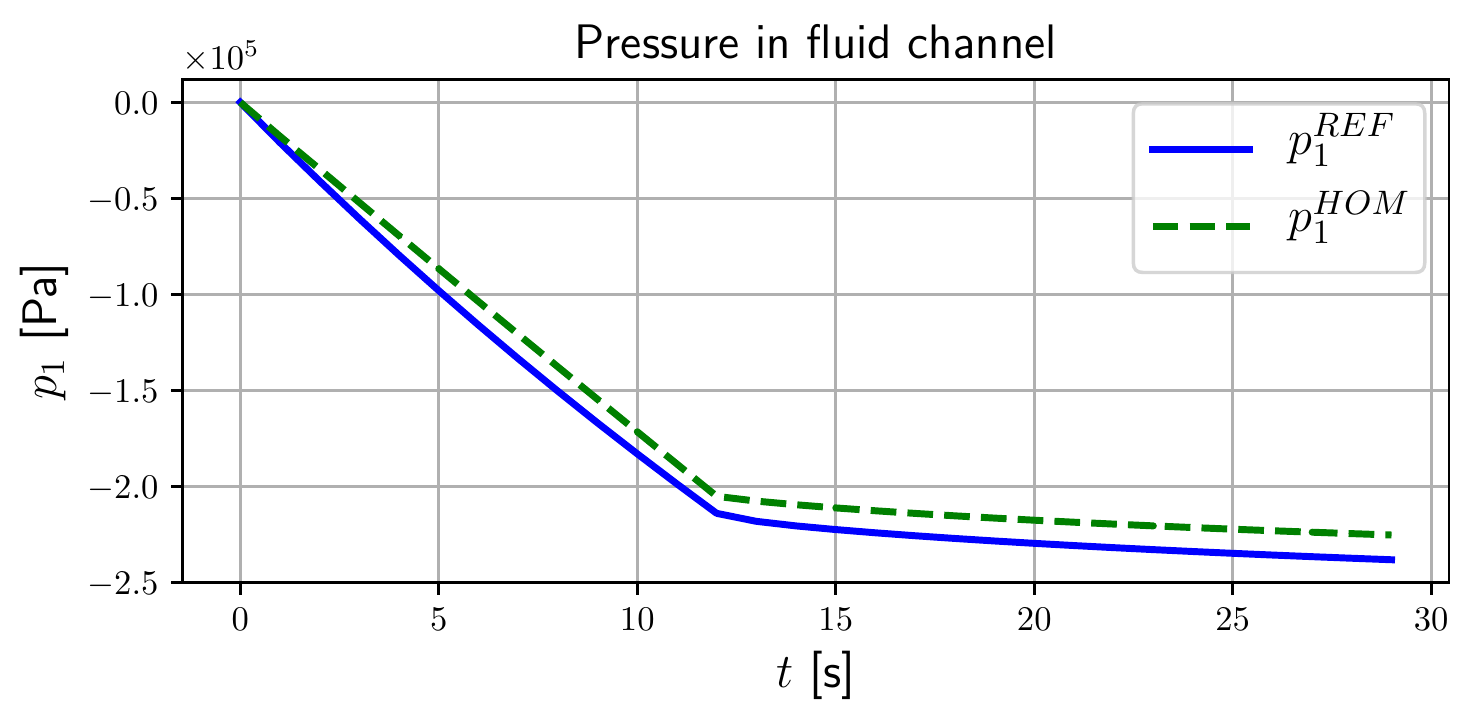}
  \caption{The pressures in the porous matrix and fluid channel obtained by
    the reference ($p^{REF}$) and homogenized ($p^{HOM}$) models.
  }\label{fig-validation_compare}
\end{figure}

The time histories of the pressures in the porous matrix and fluid channel
obtained by the reference ($p^{REF}$) and homogenized ($p^{HOM}$) models are
compared in Fig.~\ref{fig-validation_compare}. The channel pressure $p_1^{HOM}$
is the macroscopic variable involved in system \eq{eq-mac-12}, \eq{eq-mac-13},
while the matrix pressure $p_3^{HOM}$ comes from calculations on the
microscopic reference cell, see Algorithm \ref{alg1}. Note that
due to the geometry and applied boundary conditions, the computed pressure
fields are homogeneous in the matrix and channel parts.
The obtained results entitle us to conclude, that the responses of the presented
homogenized model are in a close agreement with those of the reference
model. 

\subsection{2D shear test}\label{sec:num_shear_test}

In this and the following section,
we demonstrate the proposed two-scale modelling of
large deforming porous structures using two simple illustrative simulations.

\begin{figure}
  \centering
  \begin{tabular}{p{0.5\textwidth} p{0.5\textwidth}}
    \vspace{0pt} \hfil \includegraphics[width=0.95\linewidth]{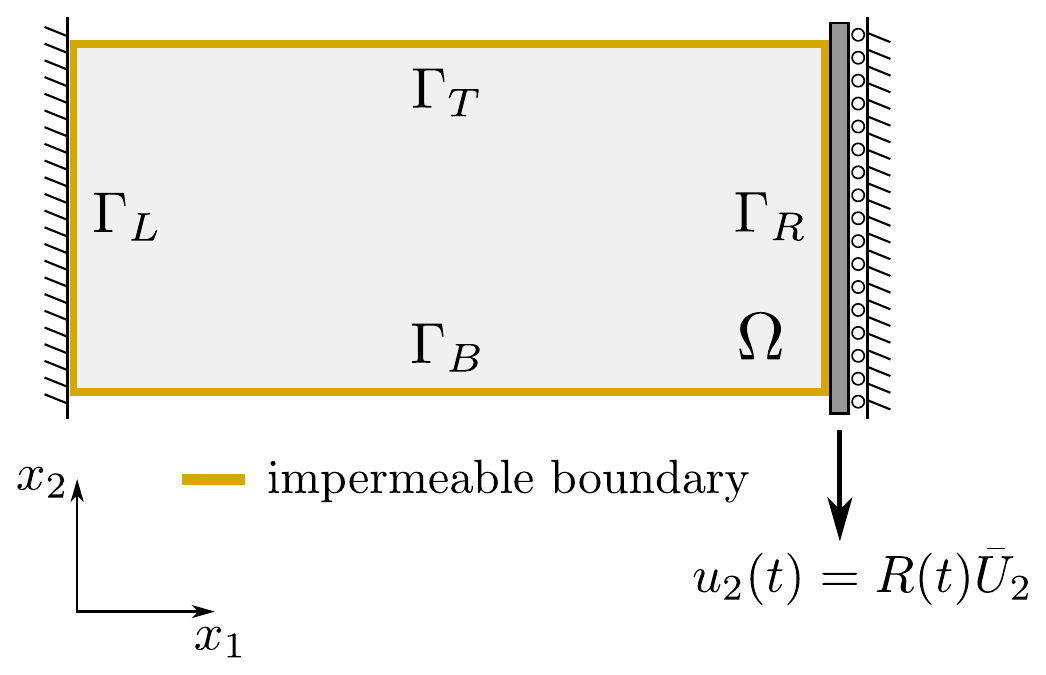}\hfil &
    \vspace{0pt} \hfil \includegraphics[width=0.6\linewidth]{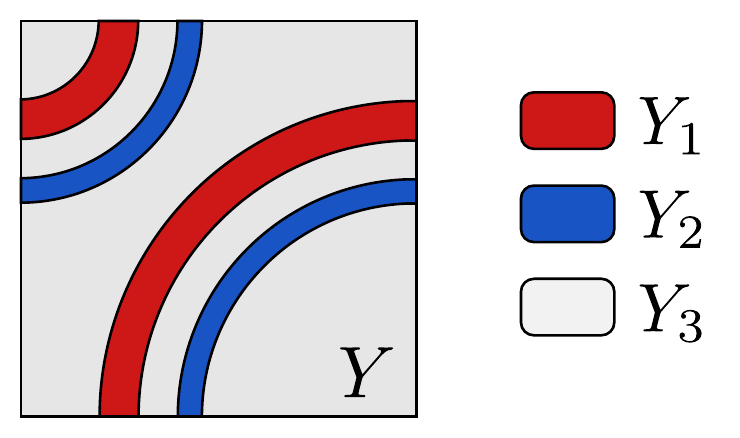}\hfil
  \end{tabular}
  \caption{Shear test: left -- boundary conditions applied to the macroscopic 2D sample;
           right -- geometry of the microscopic periodic cell.
           }\label{fig-2d_example_bc_Y}
\end{figure}

We consider the macroscopic sample of the same dimensions as in the validation
test. The applied boundary conditions, depicted in
Fig.~\ref{fig-2d_example_bc_Y} left, are as follows: the sample is fixed to the
rigid frame at the left boundary ($u_1 = u_2 = 0$ on $\Gamma_{L}$), the
displacement of the right boundary is zero in $x_1$-direction ($u_1 = 0$ on
$\Gamma_{R}$) and it is driven in time by function $u_2(t) =
R(t) \bar{U}_2$ in $x_2$-direction, where $\bar{U}_2 = -0.08$\,m and the ramp function $R(t)$ is defined
in \ref{sec:validation-test}. The reference periodic cell is established by
two curved channels embedded in the porous matrix, see
Fig.~\ref{fig-2d_example_bc_Y} right, with the material parameters defined in
Tab.~\ref{tab-mat_params}. The scaling parameter for the following simulations
is $\veps=10^{-3}$.

\begin{figure}
  \centering
  \includegraphics[width=0.99\linewidth]{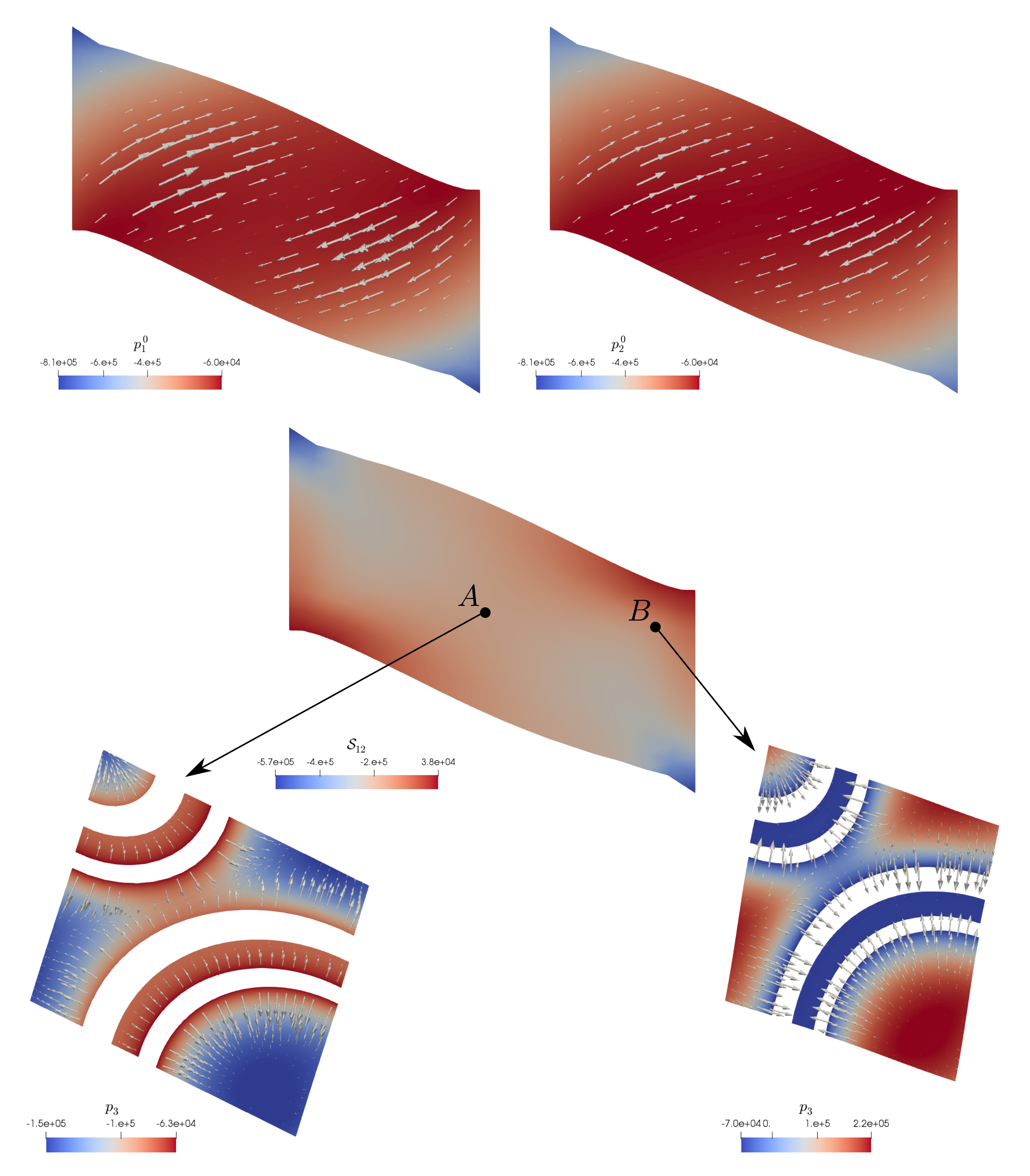}
  \caption{
    Shear test: top -- deformed macroscopic sample at time $t=0.8$\,s,
    pressure fields $p_1^0$, $p_2^0$ and the directions of perfusion velocities;
    middle -- magnitude of the averaged Cauchy stress tensor $\Scalbf$;
    bottom -- deformed microscopic reference cells at macroscopic points $A$, $B$,
    reconstructed pressure field $p_3$ associated with the solid part $Y_3$
    and the directions of perfusion velocities at the microscopic level.
  }\label{fig-2d_example_macmic}
\end{figure}

\begin{figure}
  \centering
  \includegraphics[width=0.498\linewidth]{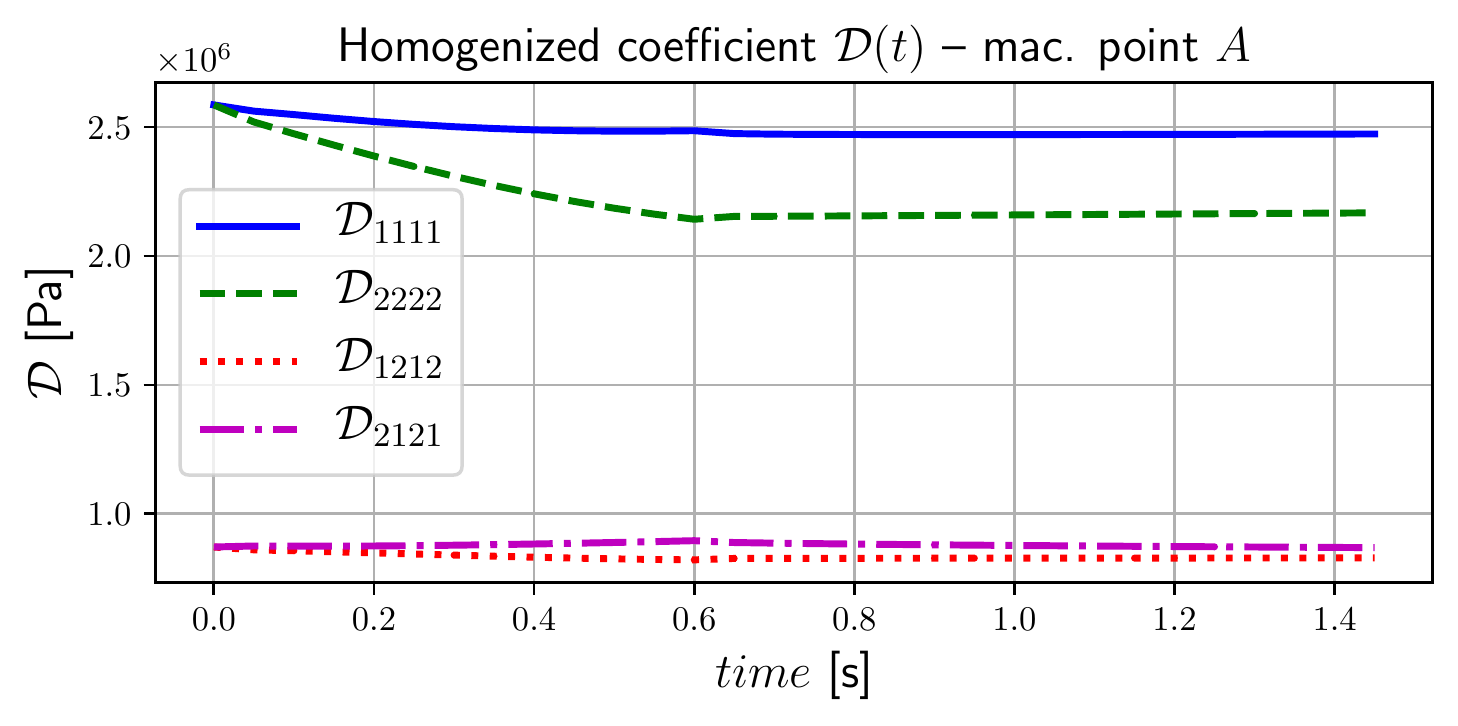}\hfil
  \includegraphics[width=0.498\linewidth]{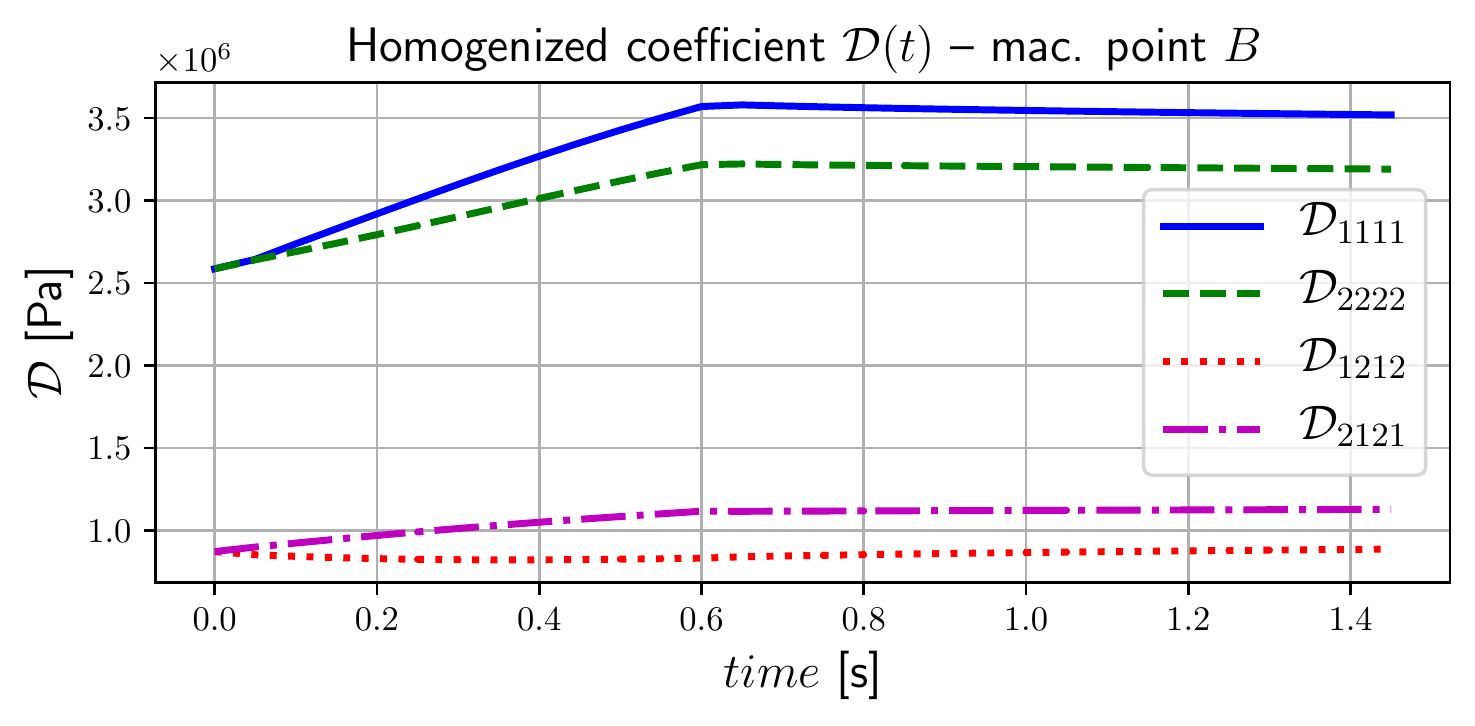}\\
  \includegraphics[width=0.498\linewidth]{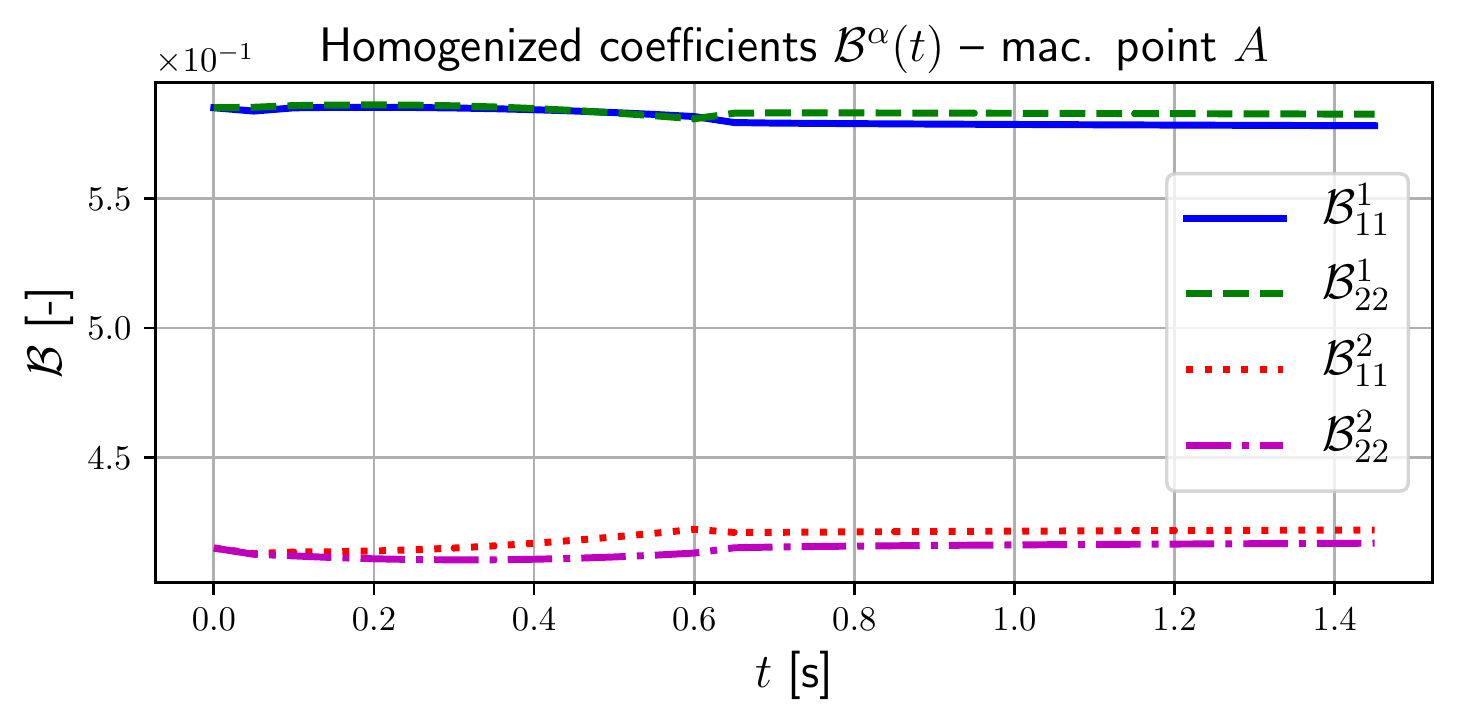}\hfil
  \includegraphics[width=0.498\linewidth]{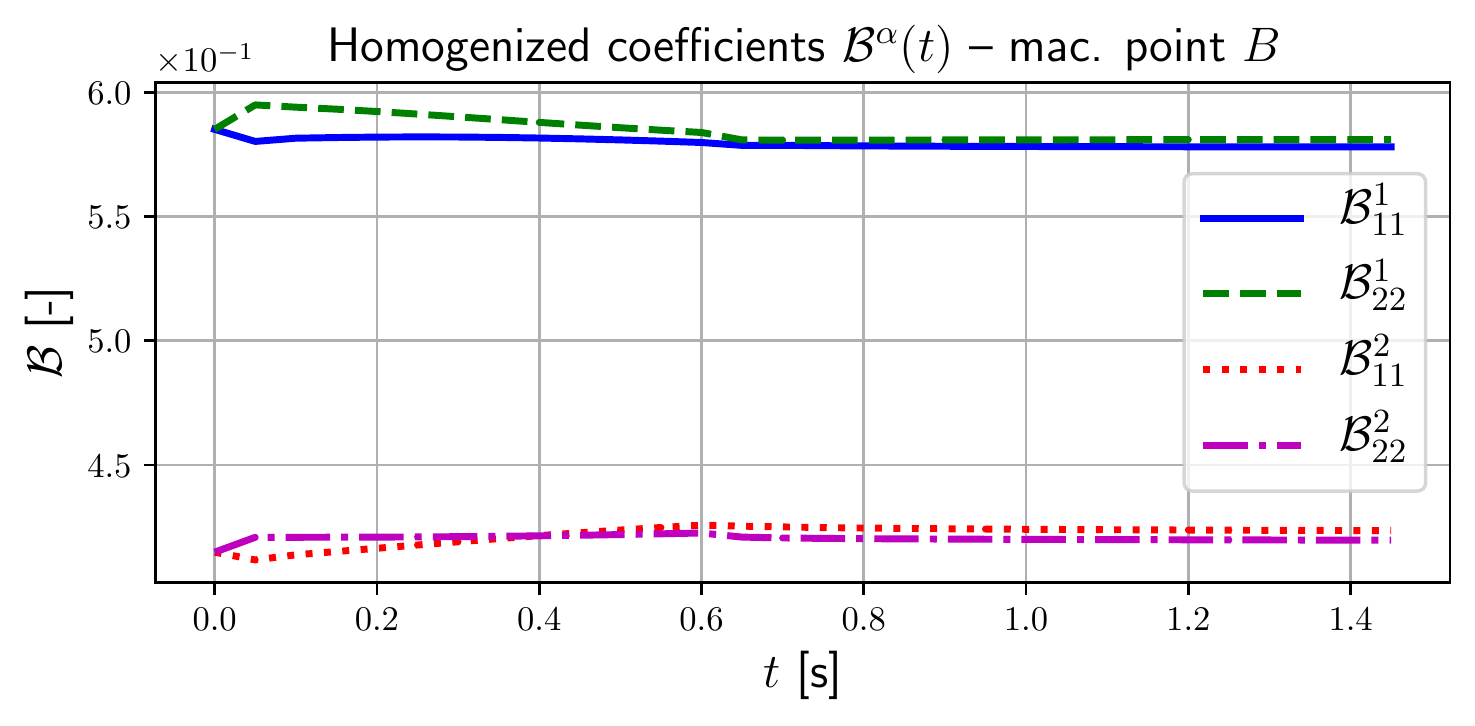}\\
  \includegraphics[width=0.498\linewidth]{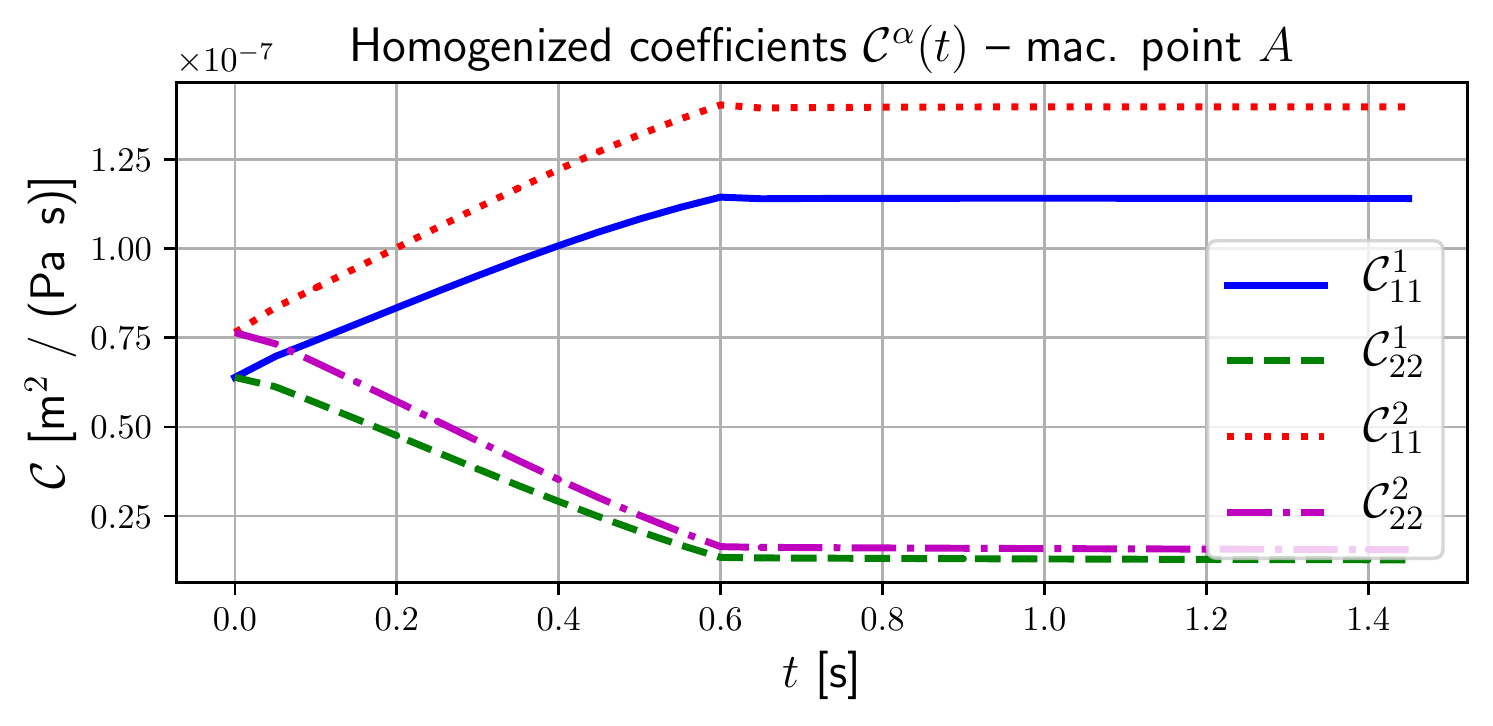}\hfil
  \includegraphics[width=0.498\linewidth]{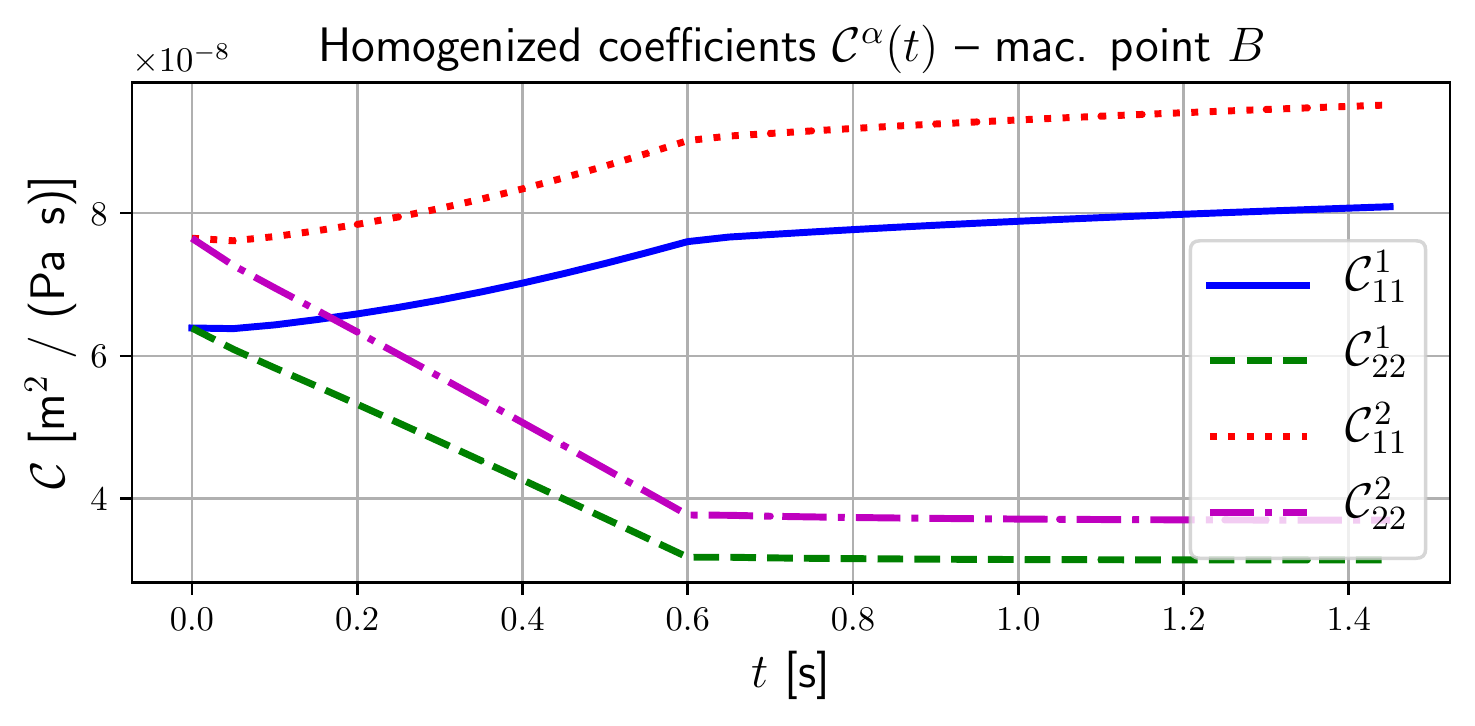}\\
  \includegraphics[width=0.498\linewidth]{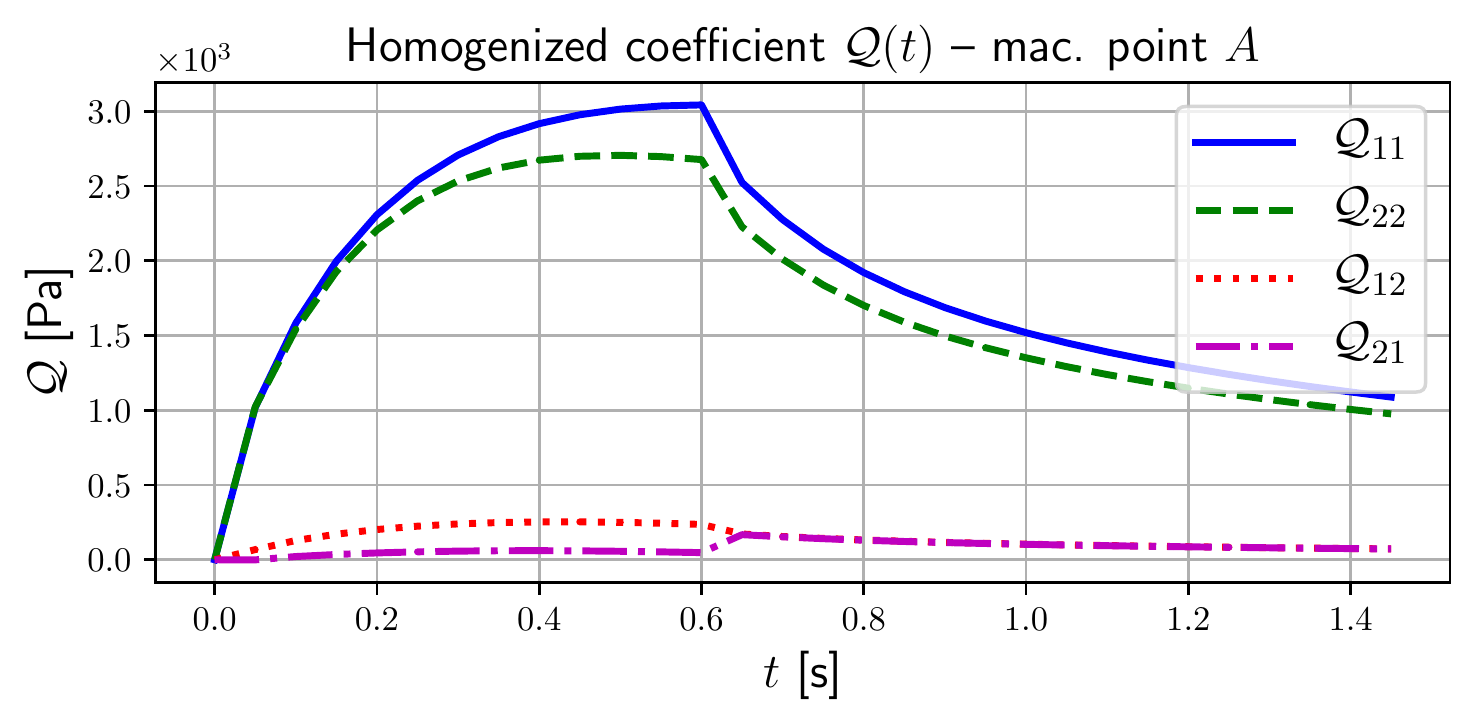}\hfil
  \includegraphics[width=0.498\linewidth]{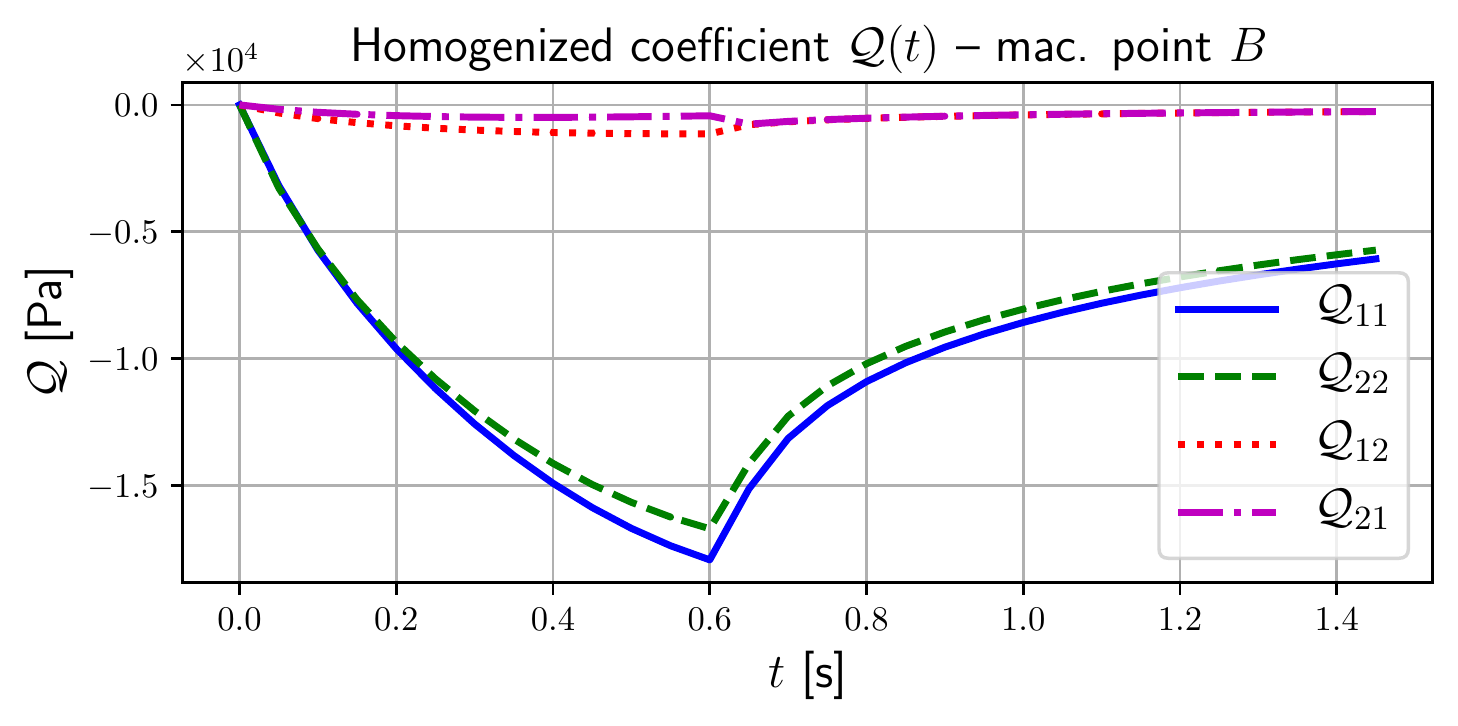}\\
  \includegraphics[width=0.498\linewidth]{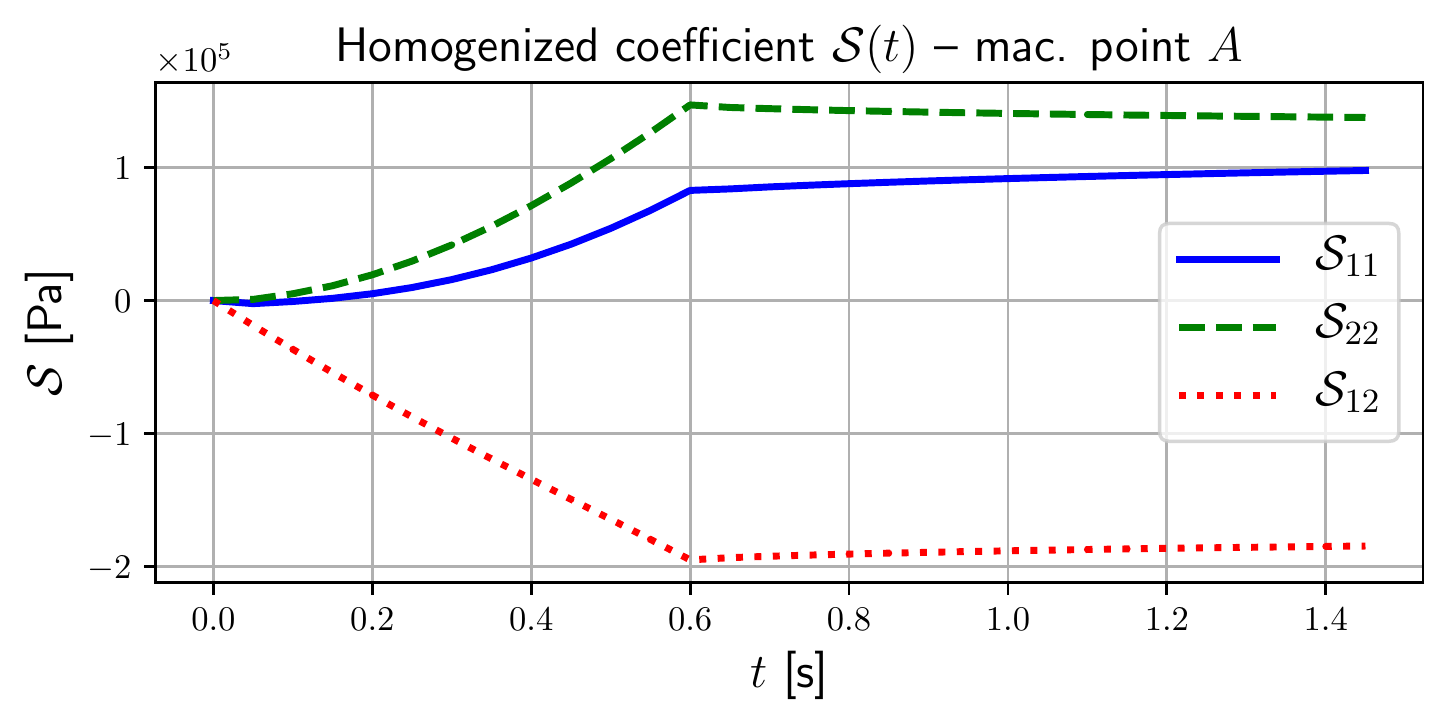}\hfil
  \includegraphics[width=0.498\linewidth]{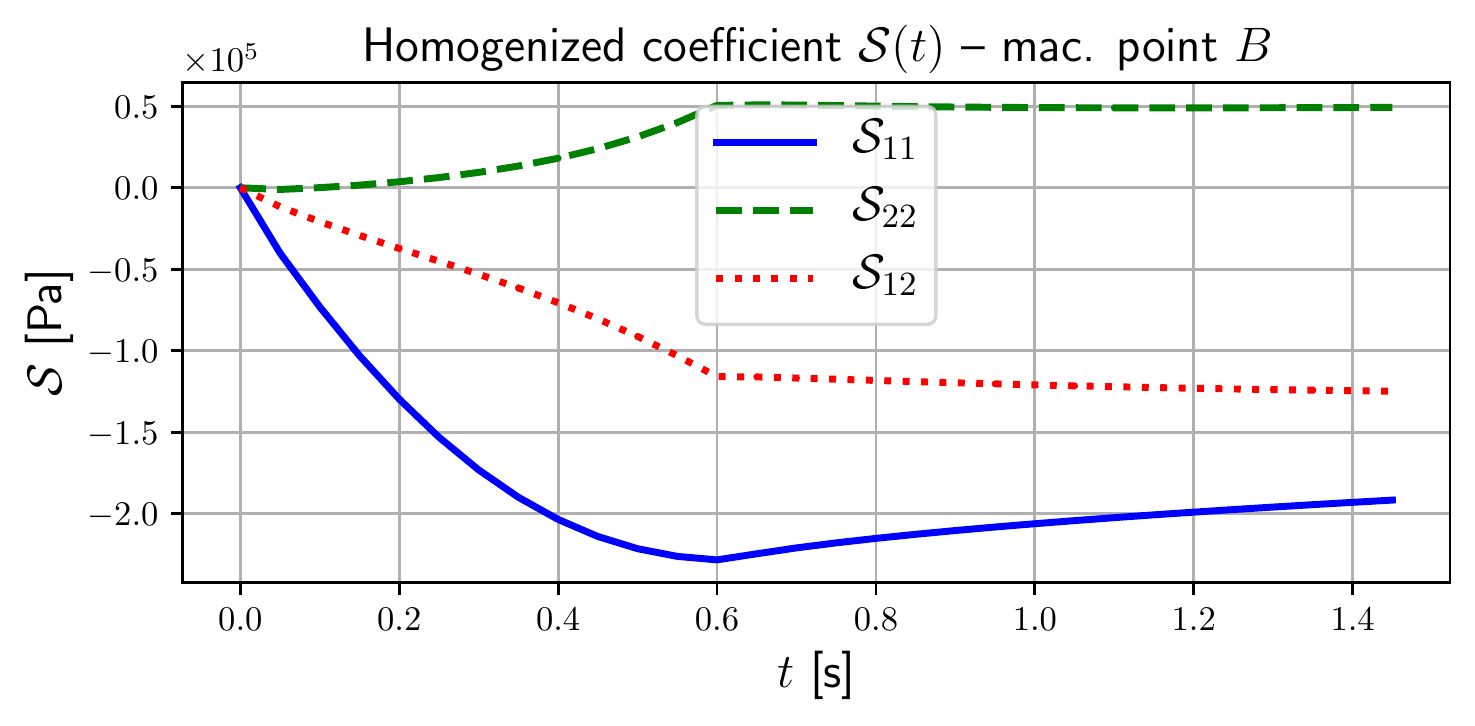}\\
  \caption{Shear test: time evolution of the homogenized coefficients $\Dcalbf$, $\Bcalbf^\alpha$,
    $\Ccalbf^\alpha$, $\Qcalbf$, $\Scalbf$ at two selected macroscopic points $A$ and $B$.}
    \label{fig-2d_example_homcf}
\end{figure}

\begin{figure}
  \centering
  \includegraphics[width=0.498\linewidth]{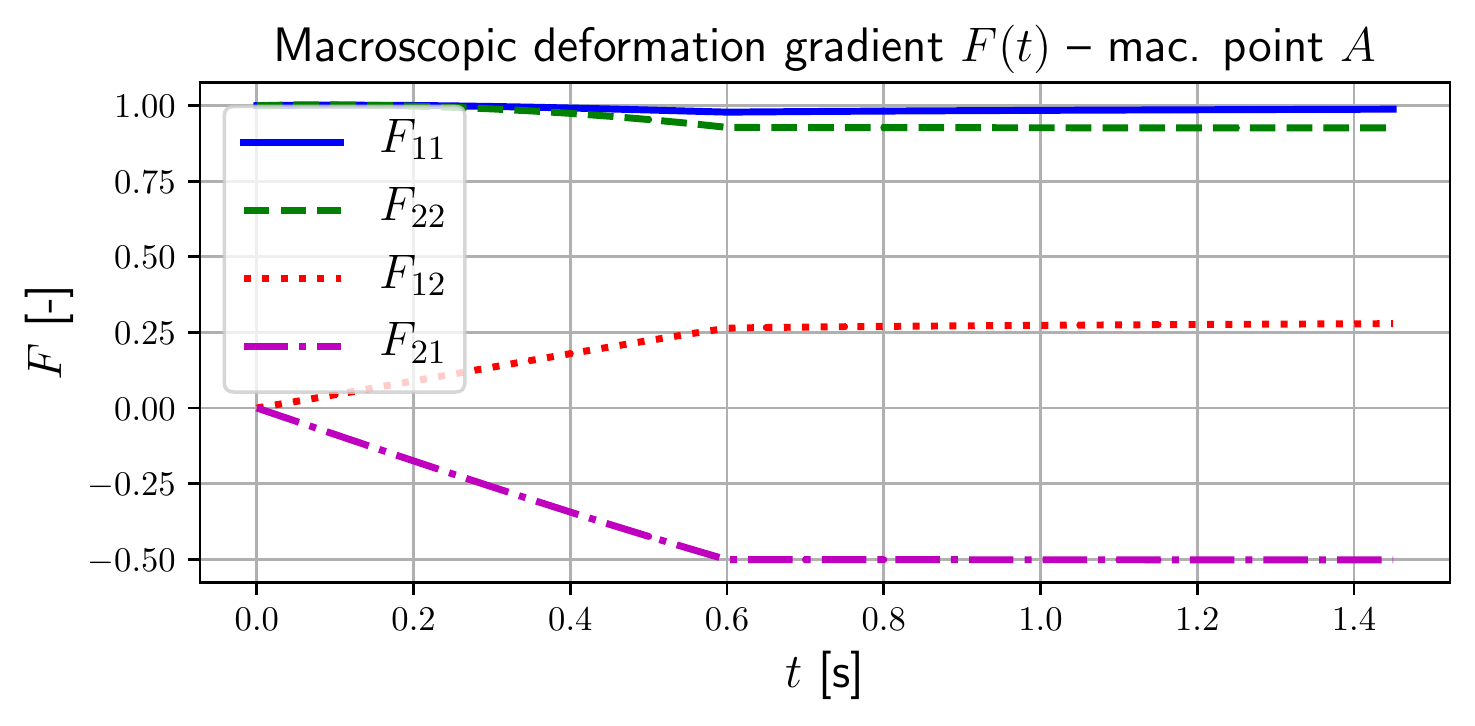}\hfil
  \includegraphics[width=0.498\linewidth]{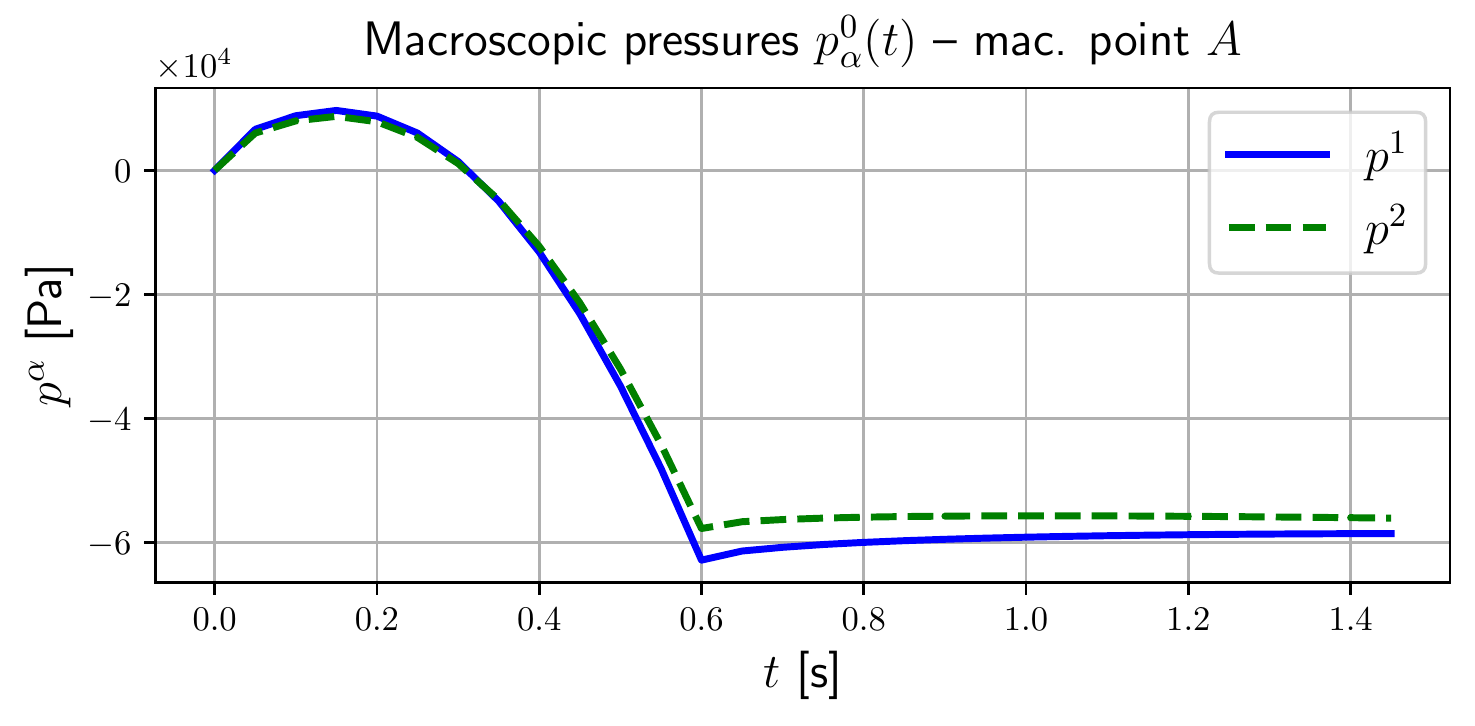}\\
  \caption{Shear test: left -- components of the macroscopic deformation gradient $\Fb$ at point A;
    right -- macroscopic pressures $p^0_\alpha$ at point A.}
    \label{fig-2d_example_def_p}
\end{figure}

The resulting macroscopic deformations, pressure distributions, averaged
Cauchy stress component $\mathcal{S}_{12}$, the deformed microscopic
cells $Y$ and the pressure field $p_3$ together with perfusion velocities $w_3$
at two selected points $A$, $B$ of domain $\Omega$ are displayed in
Fig.~\ref{fig-2d_example_macmic} for a particular time $t=0.8$\,s. In
Fig.~\ref{fig-2d_example_homcf} we show the time evolution of the homogenized
coefficients $\Dcalbf$, $\Bcalbf^\alpha$, $\Ccalbf^\alpha$, $\Qcalbf$,
$\Scalbf$ at a given macroscopic points. The components of the macroscopic
deformation gradient $\Fb$ and channel pressures $p^0_\alpha$ are shown in
Fig.~\ref{fig-2d_example_def_p}.

\subsection{2D inflation test}\label{sec:num_infl_test}

\begin{figure}
  \centering
  \begin{tabular}{p{0.38\textwidth} p{0.62\textwidth}}
    \vspace{0pt} \hfil \includegraphics[width=0.99\linewidth]{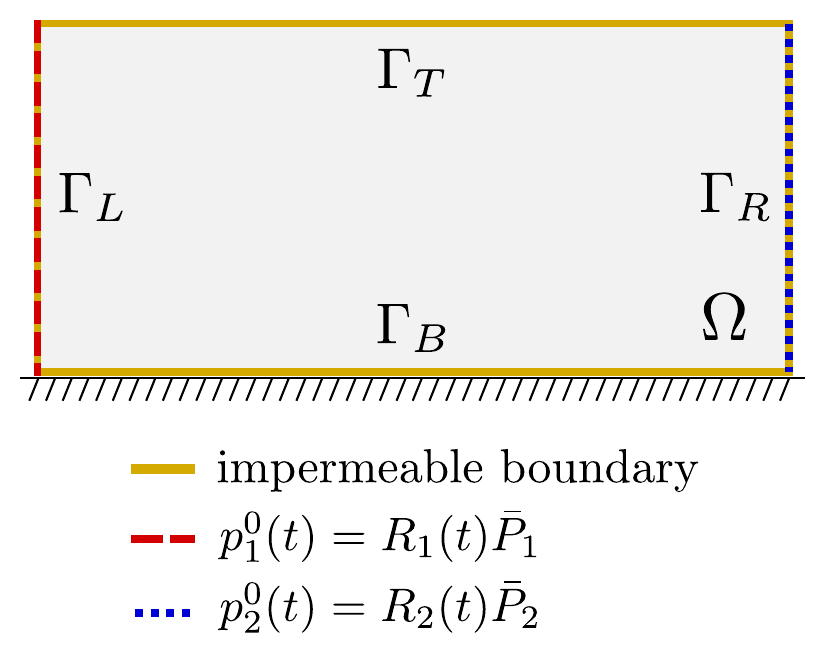}\hfil &
    \vspace{0pt} \hfil \includegraphics[width=0.99\linewidth]{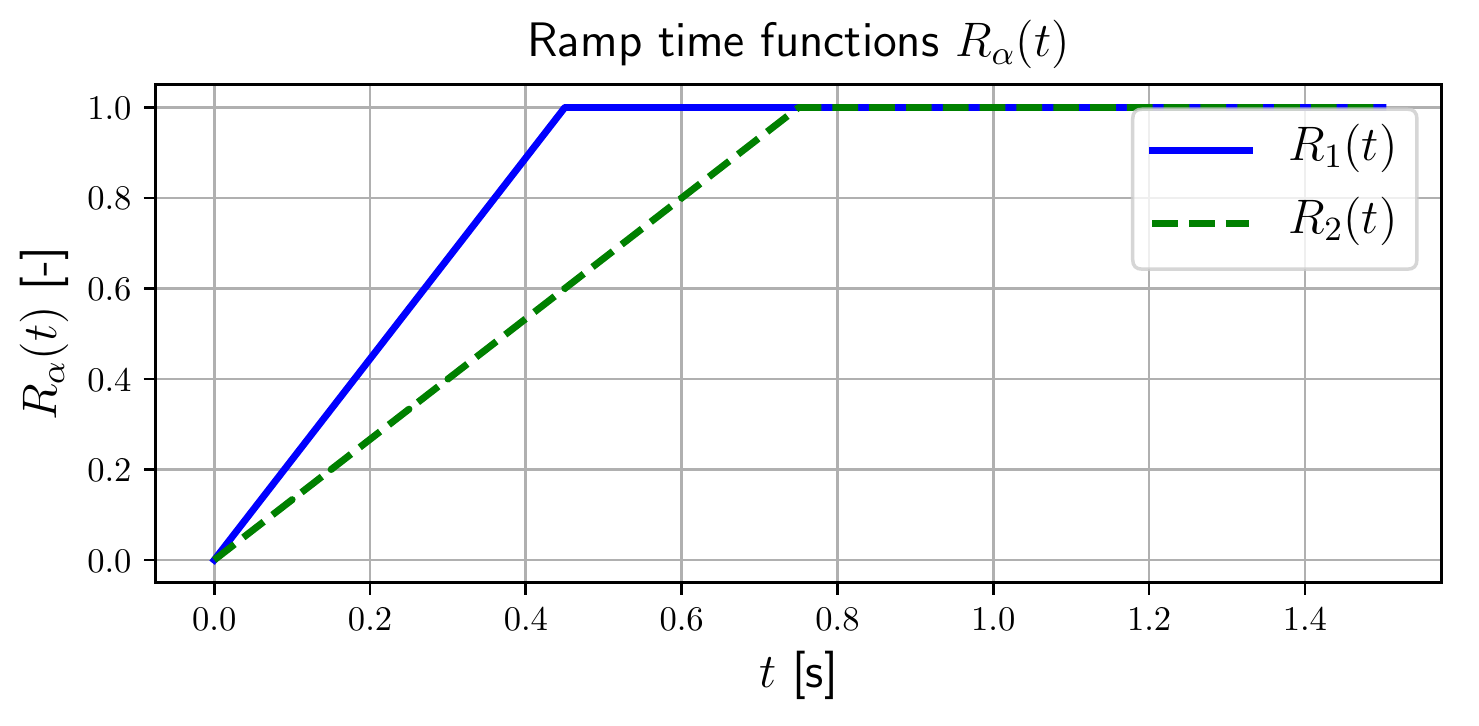}\hfil
  \end{tabular}\\
  \caption{Inflation test: left -- boundary conditions applied to the macroscopic 2D sample;
           right -- ramp functions $R_\alpha(t)$
           employed in the pressure boundary conditions.}\label{fig-2d_example2_bc}
\end{figure}

The same microscopic and macroscopic geometries and the material parameters of
the constituents as in the previous example are employed in the simple
inflation test presented in this part. The macroscopic sample is attached to
the rigid frame at the bottom edge ($u_1 = u_2 = 0$ on $\Gamma_{B}$) and two
prescribed pressures $\bar{P}_1 = 3\cdot 10^5$\,Pa and $\bar{P}_2 = 1.5\cdot 10^5$\,Pa
multiplied by the ramp functions $R_\alpha(t)$ are applied at the left
($\Gamma_L$) and right ($\Gamma_R$) edges of the sample, see
Fig.~\ref{fig-2d_example2_bc}. The part of the boundary where $p_{\alpha}$ is
not imposed, is assumed to be impermeable for fluid system $\alpha$, e.i.
$\frac{\pd p_{1}}{\pd n} =0$ on $\Gamma \setminus \Gamma_{L}$ and $\frac{\pd
p_{2}}{\pd n} =0$ on $\Gamma \setminus \Gamma_{R}$.

\begin{figure}
  \centering
  \includegraphics[width=0.99\linewidth]{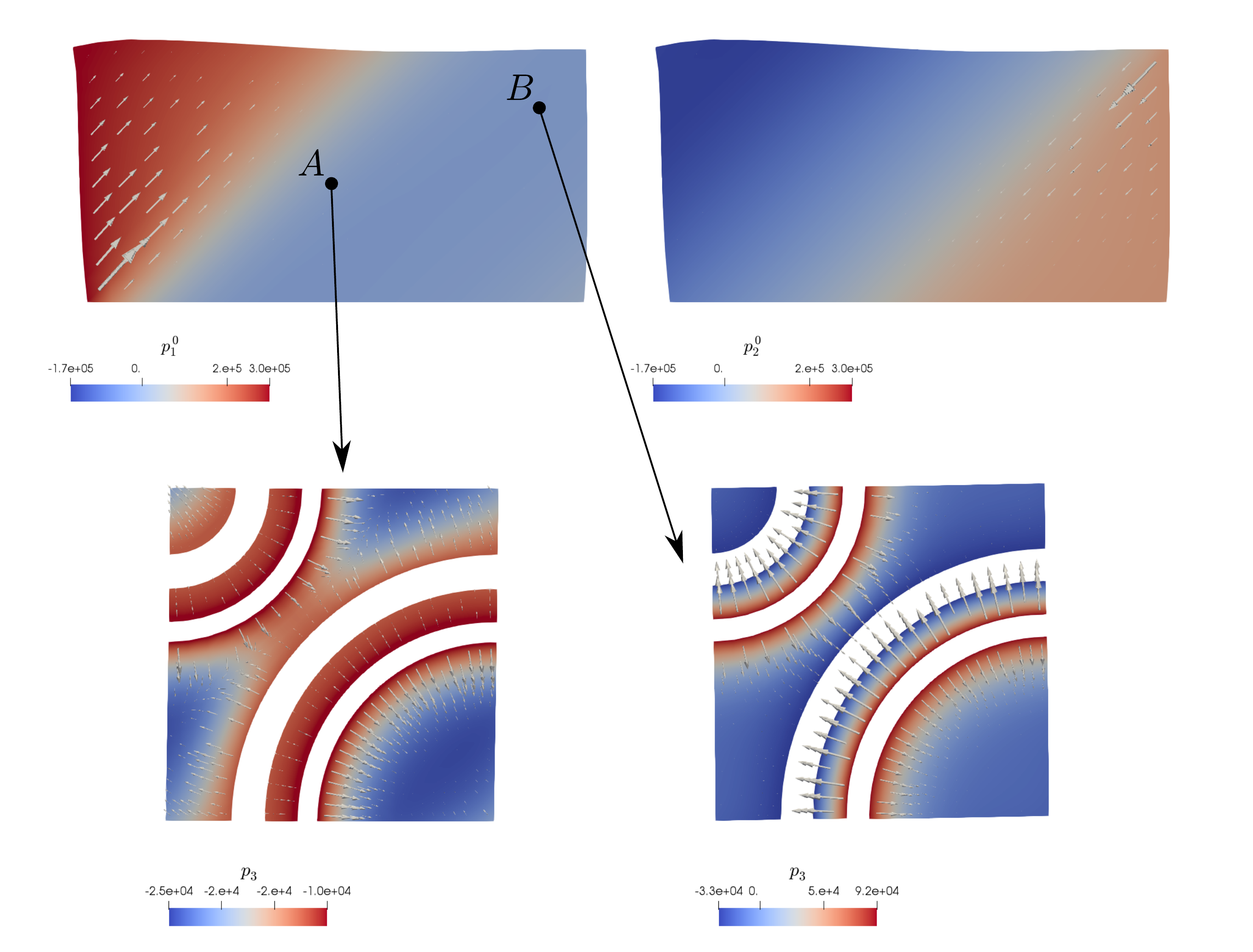}
  \caption{
    Inflation test: top -- deformed macroscopic sample at time $t=0.8$\,s,
    pressure fields $p_1^0$, $p_2^0$ and the directions of perfusion velocities;
    bottom -- deformed microscopic reference cells at macroscopic points $A$, $B$,
    reconstructed pressure field $p_3$ associated with the solid part $Y_3$
    and the directions of perfusion velocities at the microscopic level.
  }\label{fig-2d_example2_macmic}
\end{figure}

\begin{figure}
  \centering
  \includegraphics[width=0.498\linewidth]{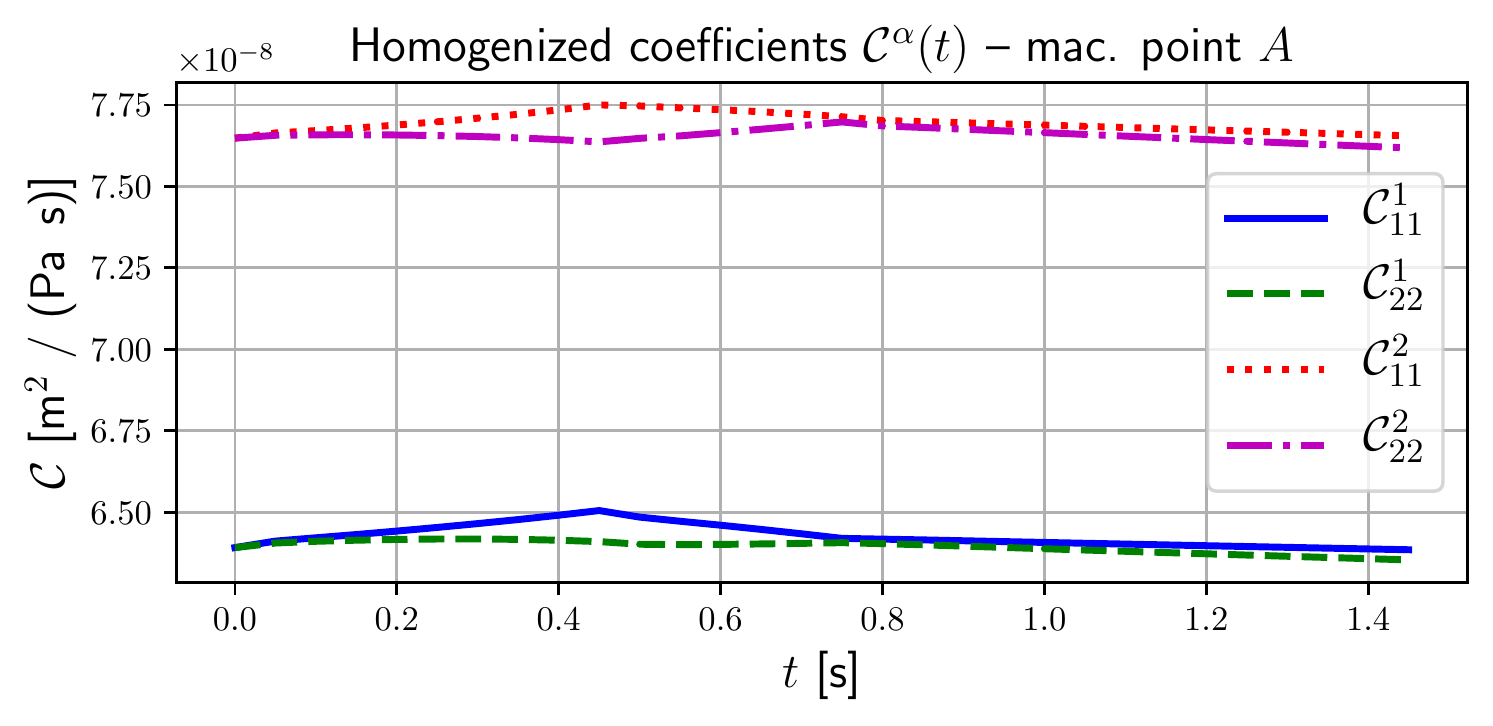}\hfil
  \includegraphics[width=0.498\linewidth]{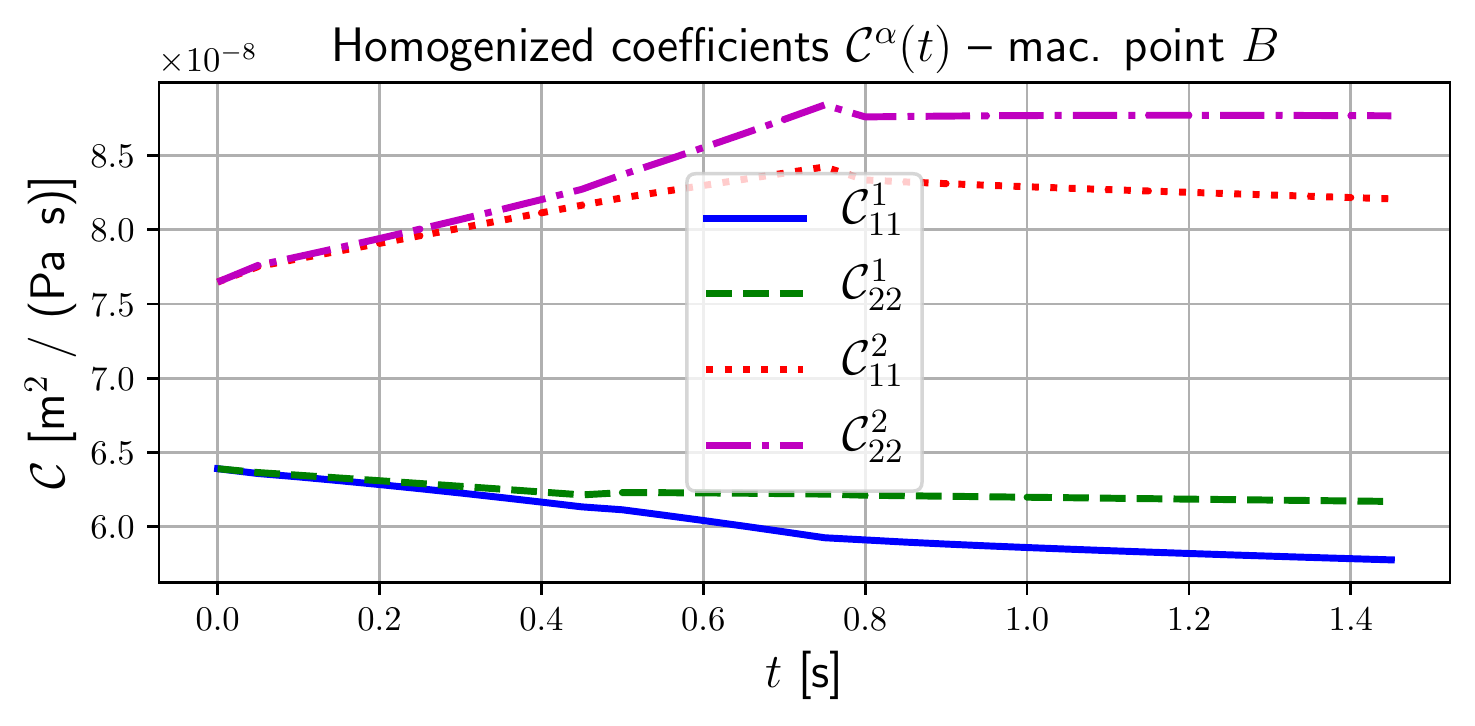}\\
  \includegraphics[width=0.498\linewidth]{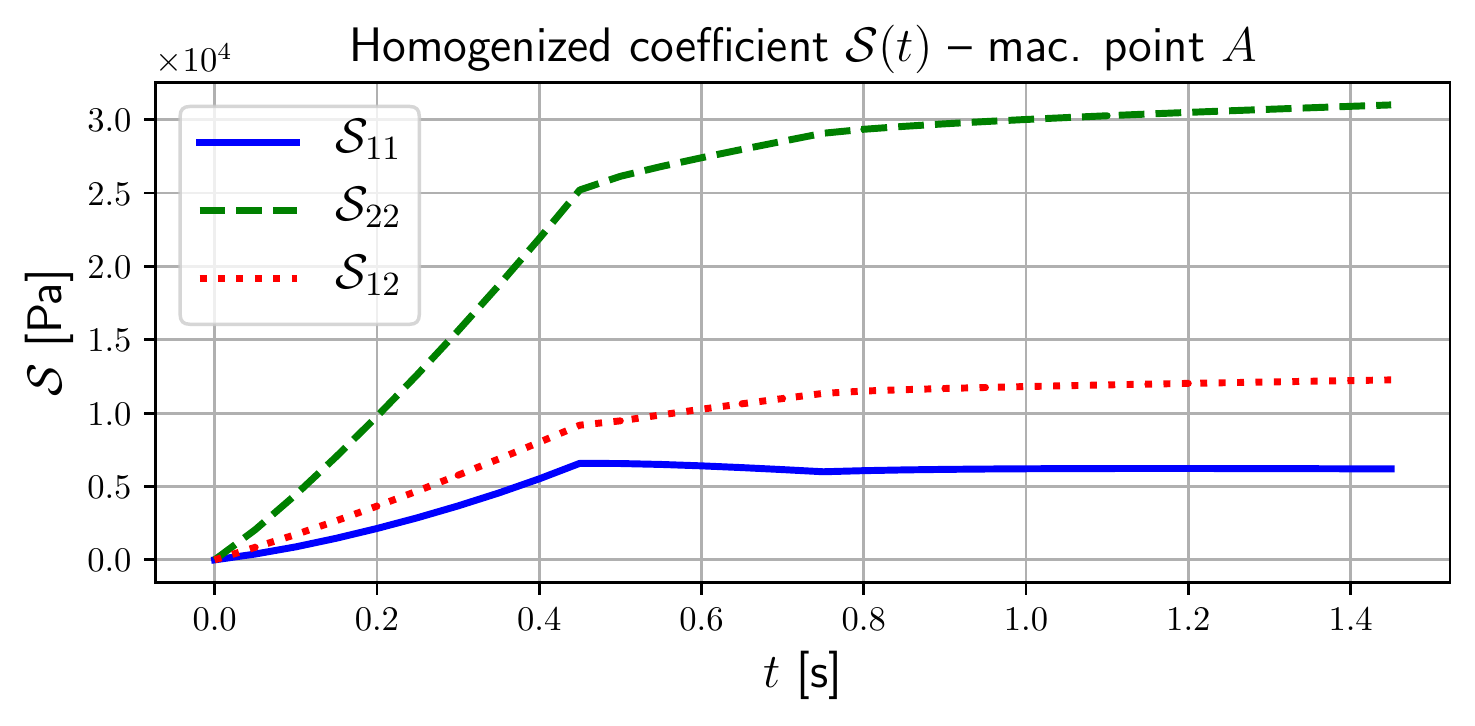}\hfil
  \includegraphics[width=0.498\linewidth]{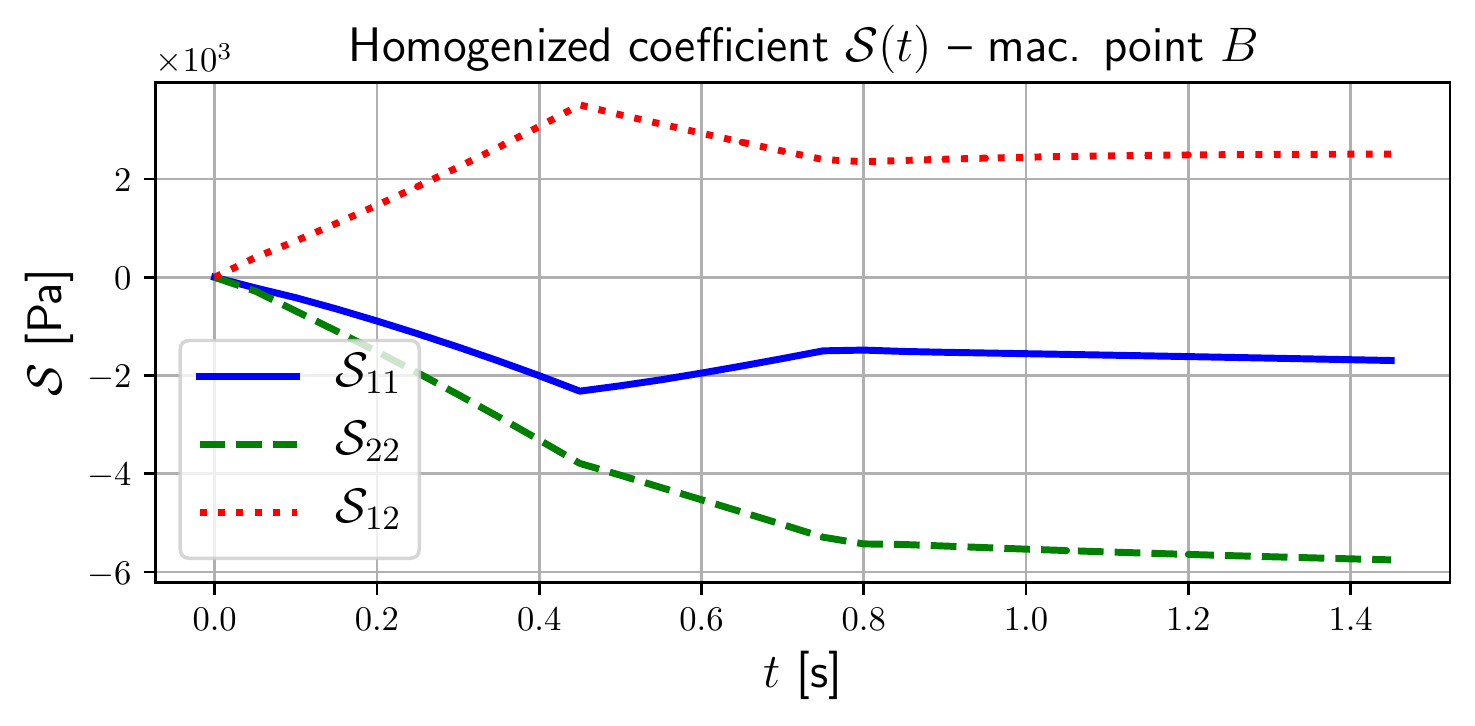}\\
  \caption{Inflation test: time evolution of the perfusion coefficients $\Ccalbf^\alpha$
   and averaged Cauchy stress $\Scalbf$ at two selected macroscopic points $A$ and $B$.}
    \label{fig-2d_example2_homcf}
\end{figure}

\begin{figure}
  \centering
  \includegraphics[width=0.498\linewidth]{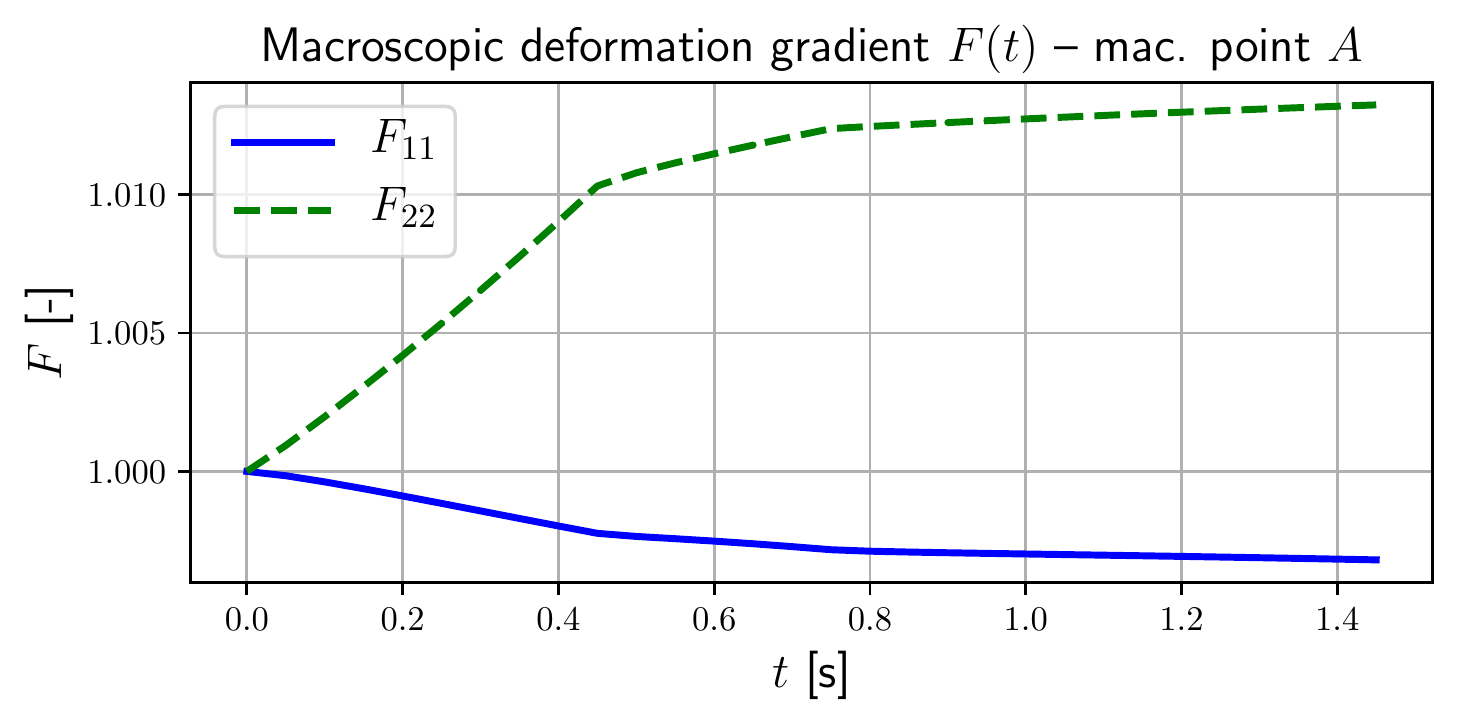}\hfil
  \includegraphics[width=0.498\linewidth]{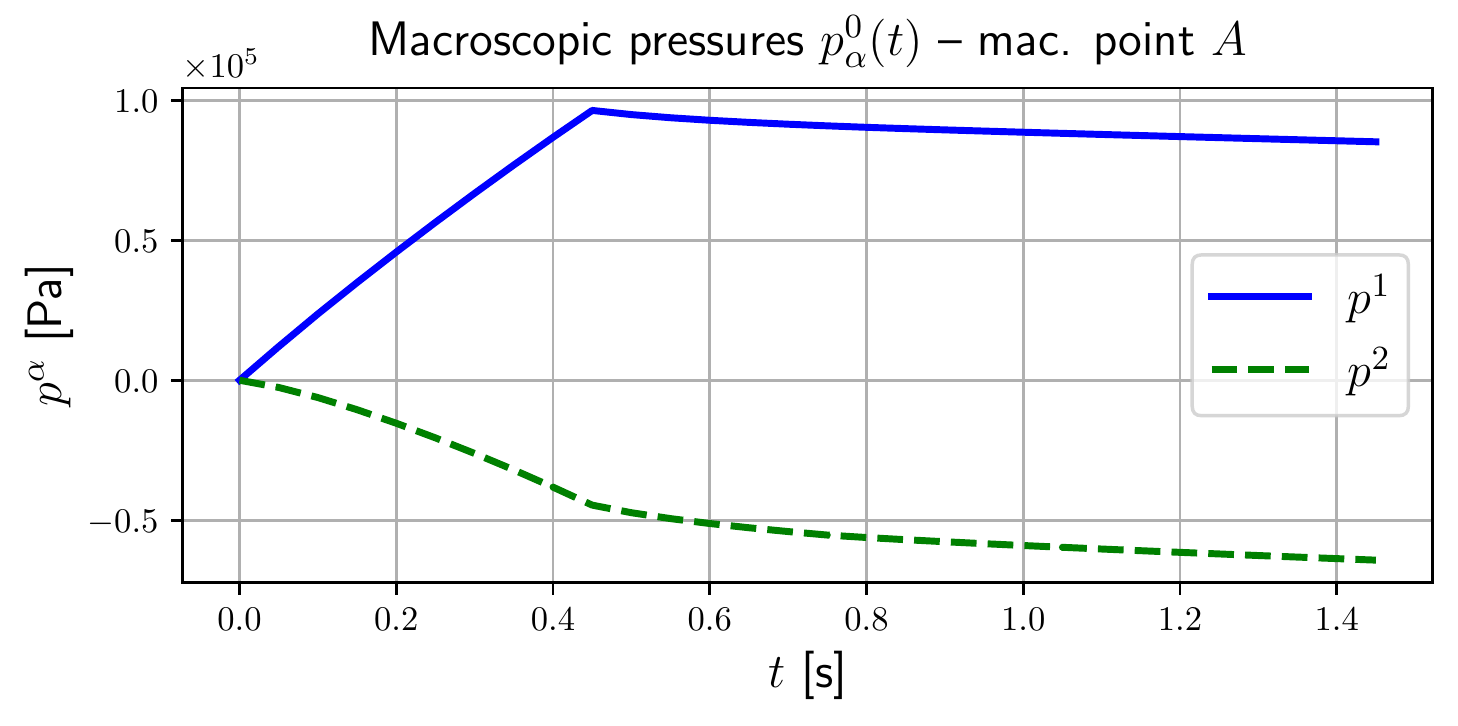}\\
  \caption{Inflation test: left -- components of the macroscopic deformation gradient $\Fb$ at point A;
    right -- macroscopic pressures $p^0_\alpha$ at point A.}
    \label{fig-2d_example2_def_p}
\end{figure}

The macroscopic and microscopic pressures and deformations induced by the
prescribed channel pressures are shown in Fig.~\ref{fig-2d_example2_macmic} for
time $t=0.8$\,s. The time changes of the homogenized coefficients $\Ccalbf$,
$\Scalbf$ is depicted in Fig.~\ref{fig-2d_example2_homcf}, deformation gradient
$\Fb$ and pressures $p^0_\alpha$ are shown in Fig.~\ref{fig-2d_example2_def_p}.

\section{Conclusion and outlook}\label{sec-conlusion}

We have proposed the multiscale model of the fluid-saturated porous media
undergoing large deformations. Its derivation is based on our previous work
\cite{Rohan-Lukes-ADES2017}, where we described modelling of large deforming
fluid saturated porous media in the Eulerian framework using a sequential
linearization concept for the Biot-type model associated with the deformed
configuration.

In this paper, we considered locally periodic heterogeneities featuring the
porous medium at the ``microscopic'' scale which, in reality, may correspond
rather to the mesoscopic scale when the multiscale, hierarchical homogenization
is in mind. In particular, we considered strong heterogeneities in the
permeability to account for the double porosity features of the medium, see \eg
\cite{Arbogast,Pride2003,rohan-cimrman-perfusionIJMCE2010,rohan-etal-jmps2012-bone,Rohan-Naili-Nguyen-ZAMM2018,rohan-etal-CaS-dp-validation},
we apply standard homogenization procedures to the linearized equations to
derive the two-scale mechanical model of large deforming porous media, assuming
a locally periodic microstructure. The proposed two-scale incremental scheme
leads to the finite element-square (FE$^2$) computational approach, where local
problems and consequently the homogenized properties are calculated by the
finite element method in a given points of a macroscopic finite element domain.
The coupled computational algorithm has been implemented in the Python based
finite element solver {\it SfePy} \cite{Cimrman-Lukes-SfePy2019}, which is
suitable for solving various multiscale and multi-physical problems. The
homogenization engine, part of {\it SfePy}, allows to evaluate all required
homogenized coefficients in an efficient way employing parallelization
possibilities of a computer system.

The presented 2D numerical examples demonstrate that the proposed two-scale
model provides a suitable computational tool for simulating fluid-saturated
porous structures with non-linear material behavior under large strains. Due to
the special topological arrangement of the microstructure involving two highly
permeable channels separated by the dual-porous matrix, this model is supposed
to be used in numerical modelling of deforming perfused tissues, such as liver
or myocardium, where the microstructural tissue arrangement and embedded vessel
networks may considerably affect the mechanical behavior at the organ level. To
be able to simulate real structures with complex 3D macroscopic geometries with
limited computational power, we need to employ some efficient strategies for
reduction of the number of local microscopic problems to be solved. In this
respect, various approaches related to this topic has been proposed, see e.g.\
\cite{Michel20036937, Fritzen2013143, Sepe2013725, Dvorak1994571, Eidel2019}.
These, however, are not directly applicable to the model featured by local
internal variables evolving with time, thus, to find a suitable model order
reduction remains an important topic for further research.


\section*{Acknowledgement}

The research has been supported by the grant project GACR 19-04956S of the Czech Science Foundation,
and in a part by the European Regional Development Fund-Project
``Application of Modern Technologies in Medicine and Industry'' {(No.~CZ.02.1.01/0.0/ 0.0/17~048/0007280)} of the Czech Ministry of Education, Youth and Sports. 

\appendix

\section{Weak formulation of deformation-diffusion problem}\label{app-weak_formulation}

The rate form of the governing equations for the deformation-diffusion problem
was derived in \cite{Rohan-Lukes-ADES2017}, where also inertia effects were
taken into account, and where a predictor--corrector time integration scheme
was proposed. Here, we present the equations without inertia terms and we consider
a simple one-step time discretization scheme.

\subsection{Rate form of the ``force--equilibrium'' equation}\label{app-weak_formulation1}

Substituting \eq{eq-wf08} and \eq{eq-wf09} into \eq{eq-wf07} yields
\begin{nalign}\label{app-weak-01}
  \incdelta\Phi_t&((\ub, p);(\vb,0)) \circ (\dot\ub, \dot p, \incdelta t \Vcal) = \\
    & \int_\Omt \Big[\nabla\Vcal(\sigmabf^\eff - p\Ib)
    + \Dop^\eff \eb(\dot\ub)
    + p\left(\nabla\Vcal + (\nabla\Vcal)^T\right)
    - (p\nabla\cdot\Vcal + \dot{p})\Ib\Big]:\nabla\vb\\
    &- \int_\Omt \left(\dot\fb \cdot \vb + \fb \cdot \vb \nabla \cdot \Vcal \right)\\
     = &\int_\Omt \Big[\Dop^\eff \eb(\dot\ub)
    + \nabla\Vcal \sigmabf^\eff
    - p\left((\nabla\cdot\Vcal)\Ib - (\nabla\Vcal)^T \right)
    - \dot p \Ib \Big]:\nabla\vb\\
   &- \int_\Omt \left(\dot\fb \cdot \vb + \fb \cdot \vb \nabla \cdot \Vcal \right).
\end{nalign}
Further, we use the time discretization \eq{eq-wf11}, the state increments \eq{eq-wf06}
and, in order to achieve the symmetry of the linearized system, we approximate
$\incdelta t \Vcal$ by the known increment $\incdelta\ub\tl{k}$ in terms
involving the unknown pressure increment $\incdelta p\tl{k+1}$
and by $\incdelta\ub\tl{k + 1}$ in the rest. Thus, we get
\begin{nalign}\label{app-weak-02}
  \incdelta\Phi_t&((\ub\tl{k}, p\tl{k});(\vb,0)) \circ (\incdelta\ub\tl{k + 1}, \incdelta p\tl{k + 1}, \incdelta t \Vcal) = \\
    &\int_\Omk \Dop^\eff \eb(\incdelta\ub\tl{k + 1})\eb(\vb)
    +  \int_\Omk\nabla (\incdelta\ub\tl{k + 1}) \sigmabf^\eff:\vb\\
    &- \int_\Omk p\tl{k}\left((\nabla\cdot\incdelta\ub\tl{k + 1})\Ib
       - (\nabla \incdelta\ub\tl{k + 1})^T \right):\nabla\vb\\
    &- \int_\Omk\Big[\incdelta p\tl{k+1}\left((\nabla\cdot\incdelta\ub\tl{k})\Ib
    - (\nabla \incdelta\ub\tl{k})^T \right)
    - \incdelta p\tl{k+1} \Ib \Big]:\nabla\vb\\
   &- \int_\Omk \left(\incdelta\fb\tl{k+1} \cdot \vb + \fb\tl{k} \cdot \vb \nabla \cdot \Vcal \right).
\end{nalign}
We add \eq{app-weak-02} and \eqm{eq-wf04}{1}, and due to the residual equation \eq{eq-wf03}, we obtain
\begin{nalign}\label{app-weak-03}
   \int_\Omk &\Big[\Dop^\eff \eb(\incdelta\ub\tl{k + 1})\eb(\vb)
    +  \nabla(\incdelta\ub\tl{k + 1}) \sigmabf^\eff:\nabla\vb
    - p\tl{k}\left((\nabla\cdot \incdelta\ub\tl{k + 1})\Ib - (\nabla \incdelta\ub\tl{k + 1})^T \right):\nabla\vb\Big]\\
    & - \int_\Omk \incdelta p\tl{k+1}\left[\left((\nabla\cdot\incdelta\ub\tl{k})\Ib - (\nabla \incdelta\ub\tl{k})^T \right)
    + \Ib\right] :\nabla\vb\\
    & = \int_\Omk \left(\fb\tl{k+1} \cdot \vb - \sigmabf\tl{k}:\nabla\vb\right).
\end{nalign}
The expression in the square brackets in the first integral of \eq{app-weak-03}
can be replaced by $\Aop \eb(\incdelta\ub\tl{k + 1})\eb(\vb)$
using the tangent elastic operator $\Aop$ introduced in \eq{eq-wf14}.
According to \eqm{eq-wf13}{1}, expression
$\left((\nabla\cdot\incdelta\ub\tl{k})\Ib - (\nabla\incdelta\ub\tl{k})^T \right)$
is equal to $\Bb(\incdelta\ub\tl{k})$ and the linearized equilibrium equation attains
the following form
\begin{nalign}\label{app-weak-04}
  \int_\Omk & \Aop \eb(\incdelta\ub\tl{k + 1})\eb(\vb)
   - \int_\Omk\incdelta p\tl{k+1}\left(\Bb(\incdelta\ub\tl{k}) + \Ib \right):\nabla\vb
   = \int_\Omk \left(\fb\tl{k+1} \cdot \vb - \sigmabf\tl{k}:\nabla\vb\right).
\end{nalign}

\subsection{Rate form of the ``fluid--content'' equation}\label{app-weak_formulation2}

The analogous steps as in \ref{app-weak_formulation1} are employed to derive
the rate form of the balance of the fluid content. In Eq.~\eq{eq-wf10}, the
time derivatives are approximated by the finite differences \eq{eq-wf11},
$\incdelta t\Vcal$ is replaced by $\incdelta\ub\tl{k}$ and the state fields
$\ub$, $p$ are rewritten employing the increments defined in \eq{eq-wf06}.
Thus, from \eq{eq-wf10} we get
\begin{nalign}\label{app-weak-05}
  \incdelta\Phi_t&((\ub\tl{k}, p\tl{k});(\boldsymbol{0}, q))  \circ (\incdelta\ub\tl{k+1}, \incdelta p\tl{k+1}, \incdelta t \Vcal) = \\
 & \int_\Omt \left(q \nabla\cdot \incdelta\ub\tl{k+1} / {\incdelta t}
  + \Kb \nabla (p\tl{k}  + \incdelta p\tl{k+1})\cdot \nabla q\right) \nabla \cdot \incdelta\ub\tl{k}\\
& + \int_\Omt \left(q \nabla\cdot (\incdelta\ub\tl{k+1} - \incdelta\ub\tl{k}) / \incdelta t
  + \Kb \nabla \incdelta p\tl{k + 1} \cdot \nabla q\right)\\
& - \int_\Omt q (\nabla \incdelta\ub\tl{k})^T: \nabla \incdelta\ub\tl{k + 1} / \incdelta t
  - \int_\Omt \Kb \left(\nabla (p\tl{k} + \incdelta p\tl{k+1}) \nabla \incdelta\ub\tl{k}\right) \cdot \nabla q\\
& - \int_\Omt \Kb \nabla (p\tl{k} + \incdelta p\tl{k+1}) \cdot \left(\nabla q \nabla \incdelta\ub\tl{k}\right) 
  + \int_\Omt \incdelta\Kb \nabla p\tl{k} \cdot \nabla q - \incdelta\Jcal\tl{k+1}(q)\\
= & 
 \int_\Omt q\left((\nabla \cdot \incdelta\ub\tl{k}) \Ib
    - (\nabla \incdelta\ub\tl{k})^T + \Ib\right): \nabla \incdelta\ub\tl{k+1} / {\incdelta t}
    - \int_\Omt q \nabla\cdot \incdelta\ub\tl{k} / \incdelta t\\
& + \int_\Omt \left(\Kb \nabla\cdot\incdelta\ub\tl{k} - (\nabla\incdelta\ub\tl{k})\Kb^T - \Kb(\nabla\incdelta\ub\tl{k})^T \right)\nabla (p\tl{k} + \incdelta p\tl{k + 1})\cdot \nabla q\\
& + \int_\Omt \Kb \nabla \incdelta p\tl{k+1} \cdot \nabla q 
+ \int_\Omt \incdelta\Kb \nabla p\tl{k} \cdot \nabla q - \incdelta\Jcal\tl{k+1}(q).
\end{nalign}
The time discretization of \eqm{eq-wf04}{2} results in
\begin{nalign}\label{app-weak-06}
  \Phi_t((\ub\tl{k}, p\tl{k});(\boldsymbol{0}, q)) = 
  \int_\Omt \left(q \nabla \cdot \incdelta\ub\tl{k} / \incdelta t
   + \Kb \nabla p\tl{k} \nabla q \right) - \Jcal\tl{k}(q).
\end{nalign}
By adding \eq{app-weak-06} to \eq{app-weak-05} and considering
$0 = \Phi_{t+\incdelta t} \approx \Phi_{t} + \incdelta\Phi_{t}$, see \eq{eq-wf03} and \eq{eq-wf05},
we obtain the linearized balance equation
\begin{nalign}\label{app-weak-07}
 \int_\Omt & q\left(\Bb(\incdelta\ub\tl{k}) + \Ib\right): \nabla \incdelta\ub\tl{k+1} / {\incdelta t}
  + \int_\Omt \left(\Hb(\incdelta\ub\tl{k}) + \Kb\right)\nabla \incdelta p\tl{k + 1}\cdot \nabla q\\
& = \Jcal\tl{k+1}(q)
 - \int_\Omt \left(\Hb(\incdelta\ub\tl{k}) + \Kb + \incdelta\Kb\right)\nabla p\tl{k}\cdot \nabla q,
\end{nalign}
where we employed $\Bb(\incdelta\ub\tl{k})$ and $\Hb(\incdelta\ub\tl{k})$ defined in \eq{eq-wf13}
and $\Jcal\tl{k+1} = \Jcal\tl{k} + \incdelta\Jcal\tl{k+1}$ contains the new flux at time level $k + 1$.

\section{}\label{app-symcf}

\subsection{Symmetric expression for coefficient \texorpdfstring{$\Dcalbf$}{D}}

Upon substituting $\vb = \omegabf^{ij}$ in \eqm{eq-mic-03}{1} and $q = \pi^{kl}$
in \eqm{eq-mic-03}{2} summation of the two equations yields the identity
\begin{nalign}
  \alin{\omegabf^{kl} + \Pibf^{kl}}{\omegabf^{ij}} + \blin{3}{\pi^{kl}}{\Pi^{ij}}
    + \dt \clin{3}{\pi^{ij}}{\pi^{kl}} = 0,
\end{nalign}
and hence we get the last equality in \eq{eq-mac-03}.  However, the minor symmetry
does not hold, $\Dcal_{ijkl} \neq \Dcal_{jikl}$, in general, as
$\omegabf^{ij}$ are not symmetric with respect to indices $i,j$, i.e.\
$\omegabf^{ij} \neq \omegabf^{ji}$.

\subsection{Equality of coupling coefficients \texorpdfstring{$\Bcalbf^\alpha$}{Ba}
 and \texorpdfstring{$\Rcalbf^\alpha$}{Ra}}

We shall employ the following identities: From
\eqm{eq-mic-03}{1} and \eqm{eq-mic-04}{2} we get subsequently
\begin{nalign}\label{eq-app-a-01}
  \alin{\omegabf^{ij}}{\omegabf^\alpha} & = -\alin{\Pibf^{ij}}{\omegabf^\alpha} + 
    \blin{3}{\pi^{ij}}{\omegabf^\alpha}\\
  & = -\alin{\Pibf^{ij}}{\omegabf^\alpha} - \dt \clin{3}{\pi^\alpha}{\pi^{ij}}.
\end{nalign}
From \eqm{eq-mic-04}{1} we have that
\begin{nalign}\label{eq-app-a-02}
  \alin{\omegabf^\alpha}{\omegabf^{ij}} & = 
  \blin{3}{\pi^\alpha}{\omegabf^{ij}} + \blin{\alpha}{1}{\omegabf^{ij}}.
\end{nalign}
Further we shall use the strong form of \eqm{eq-mic-04}{2}; by integrating by parts we obtain
\begin{nalign}\label{eq-app-a-03}
  \dt &\sumalpha\int_{\Gamma_\alpha} \left[(\tilde \Kb^3 + \tilde\Hb^3 (\bar\ub))\nabla_y \pi^{ij}\right]\cdot \nb^3 q \dSy
    - \dt \int_{Y_3} \nabla_y \cdot \left [(\tilde \Kb^3 + \tilde\Hb^3 (\bar\ub))\nabla_y \pi^{ij}\right]q \\
  & + \blin{3}{q}{\omegabf^{ij} + \Pibf^{ij}} = 0\quad \forall q \in \Hp(Y_3),
\end{nalign}
and hence
\begin{nalign}\label{eq-app-a-04}
  \nabla_y \cdot \left [(\tilde \Kb^3 + \tilde\Hb^3 (\bar\ub))\nabla_y \pi^{ij}\right]
  = (\Bb(\bar\ub)+\Ib): \nabla_y(\omegabf^{ij} + \Pibf^{ij})\quad \mbox{in } Y_3.
\end{nalign}
On ``testing'' in \eq{eq-app-a-03} by $q = \pi^\alpha$ and integrating by parts again, 
\begin{nalign}\label{eq-app-a-05}
  -\dt \int_{\Gamma_\alpha} \left [(\tilde \Kb^3 + \tilde\Hb^3 (\bar\ub))\nabla_y \pi^{ij}\right]\cdot \nb^3 \dSy
    + \dt\clin{3}{\pi^{ij}}{\pi^\alpha} 
    + \blin{3}{\pi^\alpha}{\omegabf^{ij} + \Pibf^{ij}} = 0,
\end{nalign}
and thus
\begin{nalign}\label{eq-app-a-06}
  -g_3^\alpha(\bar\ub,\pi^{ij})
    + \clin{3}{\pi^{ij}}{\pi^\alpha}
    + (\dt)^{-1}\blin{3}{\pi^\alpha}{\omegabf^{ij} + \Pibf^{ij}} = 0.
\end{nalign}
Now we can rewrite \eq{eq-mac-08}; first we use \eq{eq-app-a-02} and \eq{eq-app-a-06},
which yields
\begin{nalign}\label{eq-app-a-07}
  \Rcal_{ij}^\alpha & = |Y|^{-1} \blin{\alpha}{1}{\omegabf^{ij}+\Pibf^{ij}}
    + \incdelta t |Y|^{-1}g_3^\alpha(\bar\ub,\pi^{ij})\\
  & = |Y|^{-1} \left[ \blin{\alpha}{1}{\Pibf^{ij}}
    + \alin{\omegabf^\alpha}{\omegabf^{ij}}
    - \blin{3}{\pi^\alpha}{\omegabf^{ij}} \right.\\
  & \quad  + \dt \left. \clin{3}{\pi^{ij}}{\pi^\alpha}
    + \blin{3}{\pi^\alpha}{\omegabf^{ij} + \Pibf^{ij}} \right]\\
  & = |Y|^{-1} \left[ \blin{\alpha}{1}{\Pibf^{ij}}
    + \alin{\omegabf^\alpha}{\omegabf^{ij}} \right.\\
  & \quad \left. + \dt \clin{3}{\pi^{ij}}{\pi^\alpha}
    + \blin{3}{\pi^\alpha}{\Pibf^{ij}} \right].
\end{nalign}
Then, using \eq{eq-app-a-01}, we get
\begin{nalign}\label{eq-app-a-08}
  \Rcal_{ij}^\alpha & = |Y|^{-1}  \left[ \blin{\alpha}{1}{\Pibf^{ij}}
    + \blin{3}{\pi^\alpha}{\Pibf^{ij}}
    - \alin{\Pibf^{ij}}{\omegabf^\alpha}\right] = \Bcal_{ij}^\alpha,
\end{nalign}
which is the desired identity.


\bibliography{paper-ulf-homog-perfusion_bib}
\bibliographystyle{plain}

\end{document}